\documentclass[12pt,a4paper,oneside]{amsart}

\usepackage{amsmath,amssymb,amsthm}

\everymath{\displaystyle}

\newfam\gothfam
\font\tengoth=eufm10 \textfont\gothfam=\tengoth
\font\sevengoth=eufm7 \scriptfont\gothfam=\sevengoth
\font\fivegoth=eufm5 \scriptscriptfont\gothfam=\fivegoth

\newcommand{\F}{{\bf F}}
\newcommand{\Z}{{\bf Z}}
\newcommand{\Q}{{\bf Q}}
\newcommand{\R}{{\bf R}}
\newcommand{\Zp}{\Z_{p}}
\newcommand{\Qp}{\Q_{p}}

\newcommand{\et}{{\mbox {\tiny \'et}}}
\newcommand{\M}{\mathcal{M}}

\makeatletter\renewcommand\@biblabel[1]{#1.}\makeatother

\everymath{\displaystyle}
\def\Gal{\operatorname{Gal}}

\def\sgn{\operatorname{sgn}}

\date{\today}

        \author{J. ASSIM and A. MOVAHHEDI}
        \title{ Galois Codescent For  Motivic Tame Kernels }

\newtheorem{thm}{Theorem}[section]
        \newtheorem{proposition}[thm]{Proposition}
        \newtheorem{lem}[thm]{Lemma}
        \newtheorem{cor}[thm]{Corollary}
        \newtheorem{defi}[thm]{Definition}
        \newtheorem{rem}[thm]{Remark}

\newtheorem*{thm*}{Theorem}

        \newcommand{\fd}{\rightarrow}

        \newcommand{\coker}{\mathrm{coker} \,}
        \newcommand{\Image}{\mathrm{Im} \,}

\begin{document}

\maketitle

\date{}

 \begin{abstract}
Let    $L/F$  be a finite Galois extension of number fields with an arbitrary Galois
group  $G$.
 We give an explicit description of the kernel of the natural map on motivic tame kernels $H^2_\M(o_L, {\Z}(i))_{G} {\rightarrow}
H^2_\M(o_F, {\Z}(i))$.  
Using the link between motivic cohomology and $K$-theory, we deduce genus formulae for all even  $K$-groups 
$K_{2i-2}(o_F)$ of the ring of integers. 
As a by-product, we also obtain lower bounds
for the order of the kernel and cokernel of the functorial map
$H^2_\M(F, {\Z}(i)) \rightarrow  H^2_\M( L, {\Z}(i) )^{G}$.  
%and completely answer a queqstion raised by B. Kahn about a signature map.
\end{abstract}


%\subjclass [2000]{11R70  $K$-theory of global fields; 11R23 Iwasawa theory, 11R34 Galois cohomology,  19F27  \'Etale cohomology, higher regulators, zeta and $L$-functions.}

%\maketitle

{\bf Introduction.}  The motivic Bloch-Kato Conjecture, now a theorem of
Voevodsky, led to a more precise description of the $K$-theory of the
ring of integers of algebraic number fields.
 It  implies in particular that the Quillen-Lichtenbaum conjecture holds for odd primes. The exact deviation
 in the $2$-torsion between $K$-theory and \'etale cohomology has been determined  \cite{Kahn97, Weibel00} using   ideas in the proof of Milnor's conjecture by Voevodsky \cite{Voevodsky03}.
This leads to a better understanding of the relationship between
$K$-theory and number theory (cf e.g.  Kolster's survey
\cite{Kolster04} and  Weibel's survey \cite{Weibel05}). However,
many problems remain unsolved. Let $L/F$ be a Galois extension of
algebraic number fields with Galois group $G$. Then, for any integer
$i\geq 2$,  we are interested in the natural maps
$$  K_{2i-2} F   \rightarrow  (K_{2i-2} L)^G \; \; {\mbox {and}}  \; \;  (K_{2i-2} L ) _G \rightarrow  K_{2i-2} F.$$
By the Bloch-Kato conjecture, one can replace $K$-groups by
motivic cohomology groups and then study the kernels and cokernels
of the corresponding maps. The transfer map
$$\mathrm{tr}_i : H^2_\M( L, {\Z}(i) )_{G} \rightarrow
H^2_\M( F, {\Z}(i) $$ is an isomorphism,  up to known $2$-torsion,
so that the induced map on the rings of integers
$$\mathrm{tr}_i : H^2_\M( o_{L}, {\Z}(i) )_{G} \rightarrow
H^2_\M( o_{F}, {\Z}(i))$$ is also surjective (up to known
$2$-torsion). We work out the kernel of this map explicitly for any
Galois extension $L/F$ and obtain a genus formula  comparing the
order of the two groups $H^2_\M( o_{L}, {\Z}(i))_{G} $ and $H^2_\M(
o_{F}, {\Z}(i))$ (Theorem \ref{motivic genus}).  
%under a mild hypothesis for odd  $i$. 
The link between $K$-theory and motivic cohomology yields a genus formula for even  $K$-groups (Theorem \ref{genus K}). 
Its difficult part is a norm index as in Chevalley's genus formula for class groups, where the units are replaced by some odd  $K$-group. 
However, we manage to deal with this norm index in some special cases. For instance,  when  $L/F$  is cyclic of order  $p$,  the genus formulae for even $K$-groups involve only the ramification in the extension $L/F$ (Corollaries \ref{genus p}, \ref{genus 2 i odd}). 
The proof of the genus formulae uses a localisation sequence in motivic cohomology due to
Geisser, which relates (in our situation) the motivic cohomology
groups of  rings of integers of $F$ to the motivic cohomology of $F$
\cite[Theorem 1.1]{Geisser04} combined with  ideas used in the
\'etale setting \cite{Kolster-Movahhedi00, Kolster02,
Kolster03,Kolster-Movahhedi03,Assim-Movahhedi04, Assim-Movahhedi12}.

%Denote by  $e_v$ the ramification index of a prime  $v$ in  $L/F$. For a prime $p$, let  $T_p$  consist of the finite
%primes  $v$  of  $F$  such that  $p$  divides the ramification index  $e_v$. 
%Besides, 
Let  $\mu_p$  denotes the $p$-th roots of unity and  $E=F(\mu_p)$.
There exists a subgroup $D_F^{(i)}$ of $E^ \bullet$
(the \'etale Tate kernel  \cite{Kolster91}) such that
$$H^1_\et( F,\Z_p(i))/p   \simeq D_F^{(i)}/E^ {\bullet p}.$$
The elementary extension
$E(\sqrt[p]{D_F^{(i)}})/E$ is unramified outside $p$ and we have the following
arithmetic criterion of Galois co-descent for $H^2_\M( o_{L},
{\Z}(i))$ (which is a special case of 
Theorem \ref{motivic codescent thm}):

\begin{thm}
Let   $L/F$ be a finite Galois extension  of number fields  with Galois  group  $G$.
  Assume that $L/F$ is unramified at infinity.
Then, for $i\geq 2$ even,  the natural map tr$_i : H^2_\M(o_L,
{\Z}(i))_{G} \rightarrow H^2_\M(o_F, {\Z}(i))$ is surjective and the
following two conditions are equivalent
\begin{enumerate}
\item[(i)] the map $\mathrm{tr}_i : H^2_\M(o_L, {\Z}(i))_{G}
\rightarrow H^2_\M(o_F, {\Z}(i))$ is an isomorphism
 \item [(ii)] for every prime $p$  dividing  the ramification index $e_v$  in  $L/F$  for some finite prime  $v$, the Frobenius automorphisms  $\sigma_v(E(\sqrt[p]{D_F^{(i)}})/E)$, $v\in T_p\setminus S_p$,  are linearly independent in the $\mathbf{F}_p$-vector space  $\mathrm{Gal}(E(\sqrt[p]{D_F^{(i)}})/E)$. 
% \item [(ii)] for every prime $p$ dividing  the ramification index $e_v$  for some finite prime  $v$, the Frobenius automorphisms
%$\sigma_v(E(\sqrt[p]{D_F^{(i)}})/E)$, $v\in T_p\setminus S_p$,
%generate an $\mathbf{F}_p$-subspace of
%$\mathrm{Gal}(E(\sqrt[p]{D_F^{(i)}})/E)$ of
% dimension the cardinality of $T_p\setminus S_p$.
 \end{enumerate}
Here  $T_p\setminus S_p$  stands for the set of finite non-$p$-adic primes  $v$  of  $F$
such that  $p\mid e_v$.
\end{thm}

We then give the necessary and sufficient arithmetic conditions for the vanishing of \'etale cohomology groups  $H^2_\et(o_F[1/p], {\Z_p}(i))$  to go up along a  $p$-extension,  for a given prime  $p$.  
For  $p=2$, we are led to study the vanishing of the positive \'etale cohomology groups  $H^2_+(o_F[1/2], {\Z_2}(i))$. 
When  $i$  is odd, we have an  exact sequence  
$$0 \rightarrow  ({\Z}/2)^{\delta_i(F)} \rightarrow  H^2_+(o_F[1/2], \Z_2(i)) \rightarrow H^2(o_F[1/2], \Z_2(i)) \rightarrow 0$$
where  
$({\Z}/2)^{\delta_i(F)}$  is the cokernel of the signature map  (see \S \ref{subsection Signature}) 
$$ \mathrm{sgn}_F : H^1(F, {\Z_2}(i))/2 \rightarrow  \oplus_{v\, real}{\Z}/2.$$

In doing so, we also answer a question raised by B. Kahn in 
\cite[page 2]{Kahn97}  concerning the image of the above signature map (Theorem \ref{kahn question}). 

\begin{thm} 
%\label{kahn question} 
Let $n\geq 1$ be an integer.  Then there exists a totally real number field $F$ such that the image of  
the signature map
$\mathrm{sgn}_F$ has  $2$-rank $\rho_i=n$, for all $i\geq 2$ odd.
\end{thm} 

As we will see in Section 1, the signature map  $\mathrm{sgn}_F$  is trivial for   $i\geq 2$ even.

As far as the kernel and the cokernel of the natural restriction map
   $$f_i : H^2_\M( F, {\Z}(i) ) \rightarrow
H^2_\M( L, {\Z}(i) )^G$$
 are concerned, they are described by the $G$-cohomology of $H^1_\M( L, {\Z}(i) )$  \cite[Theorem 2.13]{Kolster04},
 which is in general difficult to compute. Using Borel's results on the abelian group structure of the odd $K$-groups
 $K_{2i-1} F$, it is not difficult to give an upper bound for the generator rank and the order of $\mathrm{coker} f_i$.
Providing lower bounds turn out to be complicated, as pointed out in the classical case of $K_2$ in  \cite{Kahn93}.
Using a cup-product argument, we construct a quotient of the Galois
cohomology group $H^2(G, H^1_\M(L, {\Z}(i)))$  (which is essentially
the cokernel of the map $f_i$) depending only on the ramification in  
$L/F$ (Theorem \ref{bounds coker})  leading to a lower bound for
$\mathrm{coker} f_i$. When  $G$ is a cyclic group, the knowledge of
the Herbrand quotient also allows us to give a lower bound for the
order of $\mathrm{ker} f_i$ (Corollary \ref{bounds ker}).

In the final section, given any integer  $i$  and any odd prime number
$p$, we provide a complete list of  $p$-extensions $L$  of  $\Q$ for
which the $p$-primary part of  $H^2_\M(o_L, \Z(i))$ (resp. $K_{2i-2} o_L)$) vanishes.  
For the prime $p=2$, we have the following (see \S 6)
\begin{thm}
Let   $L$ be a finite cyclic $2$-extension  of $\Q$ and $i\geq 2$.
Then the $2$-primary part of  $K_{2i-2}o_L$ vanishes exactly in the following cases
\begin{enumerate}
\item $L$ is totally complex and unramified outside a set of primes $\{2, \ell, \infty \}$ with
$\ell \equiv \pm 3 \pmod 8$.
\item  $L$ is totally real and 
\begin{enumerate}
\item either $2i-2\equiv 0 \pmod 8$    and  $L$  is unramified outside a set of primes $\{2,
\ell_1, \ell_2 \}$  with   $\ell_1 \not \equiv  1 \pmod 8$,  $\ell_2
\not \equiv  1 \pmod 8$ and $\ell_1 \not \equiv \ell_2 \pmod 8$.
\item or $2i-2 \equiv  4$ or $6 \pmod 8$    and  $L$  is unramified outside a set of primes $\{2,
\ell\}$  with  $\ell \equiv \pm 3 \pmod 8$.
 \end{enumerate}
 \end{enumerate}
\end{thm}

When  $2i-2 \not \equiv  0 \pmod 8$ or if  $L$  is totally imaginary, then cyclic extensions can be replaced by arbitrary Galois extensions in the above theorem. 

As far as the integers $i$  for which  $2i-2 \equiv  2 \pmod 8$  are concerned, we always have a surjective map  $K_{2i-2}o_L \to (\Z/2)^{r_1}$  for any number field  $L$. When  $L/\Q$  is a Galois  $2$-extension, then the   $2$-primary part of  $K_{2i-2}o_L$  is isomorphic to  $(\Z/2)^{r_1}$  precisely when  $L$  is unramified outside a set of primes $\{2,\ell , \infty\}$  with  $\ell \equiv \pm 3 \pmod 8$.

\section{Preliminaries}
\subsection{Notation and setting}
Notation frequently used throughout the paper will be described in
this subsection. We fix an algebraic number field  $F$ and let

%denote by  $F^{\bullet}$,  $o_F$  the multiplicative group of non-zero elements of $F$ and the ring of  integers of $F$ respectively.

\begin{tabular}{rl}
$F^{\bullet}$ & the multiplicative group of non-zero elements of $F$; \\
$o_{F}$ & the ring of  integers of $F$;\\
$r_1:=r_1(F)$ (resp. $r_2:=r_2(F)$) &the number of real (resp. complex) places of  $F$. 
\end{tabular}

We shall work with a Galois extension  $L/F$  with Galois group $G$.
For a place  $v$  of  $F$, we fix a prime  $w$ of  $L$  above $v$.
The residue fields at  $v$ and  $w$  will be denoted respectively by
$k_v$  and  $k_w$. We denote by  $e_v$  the ramification index of
$v$ in the extension $L/F$ and let $e:=\mathrm {lcm} \{e_v / \; v \; \mbox
{finite}\}$. We denote by  $q_v$  the order of  $k_v$  and for any
integer  $i \in \Z$ let
$$e_v^{(i)} :=\mathrm {gcd} (e_v, q_v^i-1).$$
If $\ell$ is the residue characteristic of the completion  $F_v$  of  $F$  at  $v$, we denote by
$e_v^\prime$ the order of the non-$\ell$-primary part of the inertia
group $I_v$  of $v$ in the extension $L/F$. Namely  $e_v^\prime$  is
the order of the factor group  $I_v/I_{v,1}$  where  $I_{v,1}$  is
the first ramification group of  $v$  in the local extension
$L_w/F_v$. In particular, when the Galois group
$G_v:=\mathrm{Gal}(L_w/F_v)$  is abelian, then
$e_v^{(i)}=e_v^\prime$ and therefore does not depend on the integer
$i$.

The notation  $S$  stands for a set of primes containing the infinite primes  $S_\infty$. 
We say that a prime $v\in S_\infty$ is ramified in $L/F$ if $v$ is a real place of $F$ which becomes complex in $L$.  
For a prime number  $p$,  we adopt the following usual notations:
%For a prime number  $p$, and a set  $S$  of primes which contains
%the infinite primes  $S_\infty$  as well as those which ramify in
%$L/F$, we adopt the following usual notations:

\begin{tabular}{rl}
$S_p$ & the set of the primes above $p$ and the infinite primes;\\
$o_{F}^{\prime}= o_{F}\left[ 1/p\right], \; o_{F}^S$ & the ring of $S_p$ (resp. $S$)-integers of $F$;\\
$U_{F}^{\prime}$ & the group of units of the ring $o_{F}^{\prime}$;\\
 ${A_F^\prime}$ (resp. ${A_F^\prime}^+$ ) & the $p$-primary part of the (resp. narrow) class group of  $o_{F}^{\prime}$;\\
$F_\infty$ & the cyclotomic  $\Z_p$-extension of  $F$;\\
 $E=F(\mu_p)$ & the field obtained by adding to our base field $F$ the  $p$-th roots of unity;\\
$F_{S_p}$ & the maximal extension of  $F$  unramified outside  $p$ and infinity;\\
$G_{S_p}(F)$ & the Galois group  $\Gal(F_{S_p}/F)$;\\
$\Delta$  &  the Galois group $\Gal(E/F)\simeq \Gal(E_\infty/F_\infty)$.
\end{tabular}

% GSp(F) be the Galois group over F of the maximal extension of F unramified outside p and infinity.

%As usual,  $F_v$  denotes the completion of  $F$  at a place  $v$  of  $F$.
%We always assume that the absolute norm  $Nv \equiv 1 \pmod p$ for all finite non-$p$-adic primes involved,  since such a condition is necessary for  $v$  to be able to ramify in a  $p$-extension of  $F$.  In particular
%$F_v = F_v(\mu_p) \simeq E_{v'}$    for such a non-$p$-adic prime  $v$  of  $F$  and any extension $v'$  of  $v$  to  $E$.

The cardinality of any finite set  $T$  will be denoted by  $\mid T
\mid$. If $n$ is a non-negative integer and $A$ is an abelian group,
we denote by $_nA$ the kernel of multiplication by $n$, and by $A/n$
the cokernel.

\subsection{Signature} 
\label{subsection Signature}

For any real place $v$  of the number field  $F$, let $i_v : F \to \mathbf{R}$ denote the
corresponding real embedding. The natural signature maps
${\mathrm{sgn}}_v : F^\bullet \rightarrow {\Z}/2$ (where
$\mathrm{sgn}_v(x)=0$  or  $1$  according to whether  $i_v(x) >0$ or not)
give rise to the following  surjective signature map
$$
\begin{array}{rrll}
%\sgn=\sgn_F : 
 F^\bullet/F^{\bullet 2} & \rightarrow & \oplus_{v\, real}{\Z}/2 \\
  x & \mapsto & ({\mathrm{sgn}}_v(x))_{v \; real} \\
\end{array}
$$
The commutative diagram 
$$
\begin{array}{cccccccccccc}
F^\bullet & \stackrel{2} {\longrightarrow} & F^\bullet & \longrightarrow &
F^\bullet/F^{\bullet 2}& \longrightarrow & 0\\
\downarrow  \oplus i_v&    &   \downarrow \oplus i_v &             &    \downarrow \\
\oplus_ {v \; \mathrm{real}} \R^\bullet & \stackrel{2} 
{\longrightarrow} &  \oplus_ {v \; \mathrm{real}} \R^\bullet& \longrightarrow &\oplus_ {v \; \mathrm{real}} \R^\bullet /\R^{\bullet 2}& \longrightarrow &  0 , 
  \end{array}
  $$
together with Kummer theory, shows that the above signature map is the same as the localisation map 
$$H^1(F, \mu_2) \longrightarrow  \oplus_ {v \; \mathrm{real}} H^1(F_v, \mu_2),$$ 
where  $H^r(K, \quad)$  denotes the (continuous) Galois cohomology of the absolute Galois group of any field $K$. 

We will also be considering restrictions of the signature map such as the following composite map  
$$ H^1(F, {\Z_2}(i))/2 \hookrightarrow H^1(F, {\Z}/2(i)) \simeq F^\bullet/F^{\bullet 2}  \rightarrow  \oplus_{v\, real}{\Z}/2  
$$
and denote by  $\delta_i := \delta_i(F)$  the rank of its cokernel: 
\begin{equation}\label{sgn}
H^1(F, {\Z_2}(i))/2 \stackrel{\mathrm{sgn}_F} \longrightarrow ({\Z}/2)^{r_1} \longrightarrow ({\Z}/2)^{\delta_i} \longrightarrow 0. 
\end{equation}

For a   finite extension $L/F$   of number fields, let $R(L/F)$ be
the set of infinite primes of $F$  which ramify in $L$  and
$r:=r(L/F):= \mid R(L/F) \mid$. Then we have a (partial) surjective
signature map  
$$
\begin{array}{cccc}
   F^\bullet/F^{\bullet 2} & \rightarrow & \oplus_{v \in R(L/F)} {\Z}/2 \\
   x & \mapsto & ({\mathrm{sgn}}_v(x))_{v \in R(L/F)} \\
\end{array}
$$
and we denote by  $s_i: =s_i(L/F)$  the rank of  the cokernel of its restriction  $\mathrm{sgn}_{L/F}$  to
$H^1(F, {\Z_2}(i))/2$ (\cite{Kolster02, Kolster03, Kolster04}): 
\begin{equation}\label{sgn s_i}
H^1(F, {\Z_2}(i))/2 \stackrel{\mathrm{sgn}_{L/F}} \longrightarrow ({\Z}/2)^{r(L/F)} \longrightarrow ({\Z}/2)^{s_i(L/F)} \longrightarrow 0. 
\end{equation}

In general, we have  $\delta_i \geq s_i$ as it can readily be seen by the following commutative diagram  
$$ \begin{array}{cccccccccccc}
H^1(F, {\Z_2}(i))/2& \stackrel{\mathrm{sgn}_F} \longrightarrow & ({\Z}/2)^{r_1} & \longrightarrow & ({\Z}/2)^{\delta_i}& \longrightarrow & 0\\
\downarrow   &    &   \downarrow  &             &    \downarrow \\
H^1(F, {\Z_2}(i))/2& \stackrel{\mathrm{sgn}_{L/F}} \longrightarrow & ({\Z}/2)^r & \longrightarrow & ({\Z}/2)^{s_i}& \longrightarrow & 0  
\end{array} $$
where the right vertical map, induced by the two others, is surjective. 
More precisely, we have   $r_1 - r \geq  \delta_i - s_i \geq 0$. 

Suppose  now that  $i$  is even and consider the following commutative diagram 
\begin{equation}
\label{diagram sgn}
\begin{array}{cccccccccccc}
0 & \rightarrow & H^1(F, {\Z_2}(i))/2& \rightarrow & H^1(F, \Z/2(i)) & \rightarrow &
{_2H^2(F, {\Z_2}(i))}& \rightarrow & 0\\
 &  & \downarrow  &    &   \downarrow  &             &    \downarrow \\
 &  & 0 & \rightarrow & \oplus_ {v \; \mathrm{real}} H^1(F_v, \Z/2(i)) & \stackrel{\cong} 
{\rightarrow} &  \oplus_ {v \; \mathrm{real}} \; {_2H^2(F_v, {\Z_2}(i))}& &  
  \end{array}
\end{equation}
where the vertical maps are the localisation maps and the isomorphism on the bottom row comes from the following two facts for any real prime  $v$ and any even $i$ \cite[page 231]{Kolster04}:\\
(i) the vanishing of  $H^1(F_v, {\Z_2}(i))$; \\
(ii) the cohomology groups  $H^1(F_v, \Z/2(i))$ and  ${_2H^2(F_v, {\Z_2}(i))}=H^2(F_v, {\Z_2}(i))$ are both cyclic of order $2$. \\
As noticed before, the middle vertical map is the signature map so that 
\begin{proposition} \label{sgn triviality}
(compare to \cite[Lemma 2.5]{Kolster03})
For any number field  $F$ and any even integer  $i$, the signature map 
$$\mathrm{sgn}_F : H^1(F, {\Z_2}(i))/2 \longrightarrow  ({\Z}/2)^{r_1}$$ 
is trivial. 
\qed
\end{proposition} 

The above proposition shows in particular that  $s_i=r$ in the exact sequence  (\ref{sgn s_i})  whenever $i$  is even. 
For  $i \equiv 3 \pmod 4$, a question \cite[page 2]{Kahn97} raised by B. Kahn concerns the $2$-rank  $\rho_i(F)$  of the image of  $\mathrm{sgn}_F$. More precisely, it is asked whether  $\rho_i(F):=r_1(F)-\delta_i = 1$  when  $r_1(F) \geq 1$.  
In other words, is the image of  $H^1(F, {\Z_2}(i))$  in  $F^\bullet/F^{\bullet 2}$  contained in  
$\{-1, +1\} \times \mbox {\{totally positive elements\}}$?

As we will see in Section \ref{S galois codescent}, for any integer  $n \geq 1$ and any odd integer  $i \geq 3$, there exists a real number field  $F$  such that   
$\rho_i(F)=n$.  

When  $i$  is odd we will often make the following hypothesis relative to the extension $L/F$ 
$$(\mathcal{H}_i) \qquad \qquad  \hat{H}^n(G, ({\Z}/2)^{\delta_i(L)})=0   \quad \mbox{ for } n = -2  \mbox{ and } -1.$$ 

This hypothesis obviously holds in the following cases: 
%(i) the number field  $L$  is totally imaginary or more generally when the signature map  $\mathrm{sgn}_L : H^1(L, {\Z_2}(i))/2 \longrightarrow  ({\Z}/2)^{r_1(L)}$  is trivial, 
%(ii)  when the signature map  $\mathrm{sgn}_L : H^1(L, {\Z_2}(i))/2 \longrightarrow  ({\Z}/2)^{r_1(L)}$  is trivial, 
%
%The above proposition shows in particular that  $s_i=r$ in the exact sequence  (\ref{sgn s_i})  whenever $i$  is even. 
%When  $i$  is odd we will make the following hypothesis "relative to"  the extension $L/K$ 
%$$(\mathcal{H}) \qquad \qquad   \hat{H}^n(G, ({\Z}/2)^{\delta_i(L)})=0   \mbox{ for all } n \in \Z.$$ 
%This hypothesis obviously holds in the following cases: \\
%(i) the number field  $L$  is totally imaginary or more generally when the signature map  
%$\mathrm{sgn}_L : H^1(L, {\Z_2}(i))/2 \longrightarrow  ({\Z}/2)^{r_1(L)}$  is trivial, \\
%(ii)  when the signature map  $\mathrm{sgn}_L$   is surjective. 
\begin{enumerate}
\item[(i)] the number field  $L$  is totally imaginary or more generally the signature map  $\mathrm{sgn}_L : H^1(L, {\Z_2}(i))/2 \longrightarrow  ({\Z}/2)^{r_1(L)}$  is trivial, 
\item[(ii)]  the signature map  $\mathrm{sgn}_L$   is surjective. 
\end{enumerate}

We will also need to consider the restriction of the signature map to  $S$-units  $U_S$  for any set  $S$  of primes of  $F$  containing the infinite primes 
$$ U_S/U_S^2 \longrightarrow  ({\Z}/2)^{r_1}.$$ 
Let  $A_S(F)$ (resp.  $A_S(F)^+$) be the  $2$-primary part of the (resp. narrow) $S$-class group of  $F$.
We then have the following well-known exact sequence 
 \begin{equation}\label{5-term}
0\rightarrow U_S^+/U_S^2 \rightarrow  U_S/U_S^2 \rightarrow  ({\Z}/2)^{r_1} 
\rightarrow  A_S(F)^+\rightarrow  A_S(F) \rightarrow 0
\end{equation}
where  $U_S^+$ is the subset of totally positive $S$-units. 
In particular,  $A_S(F) \cong A_S(F)^+$ precisely when the restriction of the signature map to   $S$-units is surjective.

\subsection{Motivic cohomology and $K$-theory}
Algebraic $K$-theory and motivic cohomology can be thought of as
global cohomology theories for  \'{e}tale cohomology.
%There have been several definitions of motivic cohomology groups.
For a survey on motivic cohomology, algebraic $K$-Theory and  their
connection with number theory we refer the reader to
\cite{Kolster04, Weibel05}. Let  $i$  be a fixed integer. The
motivic cohomology groups $H^k_\M (X, {\Z}(i))$ for a smooth scheme
$X$ over a base $B$ can be defined as the hypercohomology of Bloch's
cycles $ {\Z}(i)$ for the Zariski topology. We will be dealing with
the case where  $X=B$ is the  spectrum  of a field $K$  or of the
ring of $S$-integers   $o_F^S$  of a number  field $F$.
 We denote  the corresponding motivic  cohomology groups respectively by
$H^k_{\M} (K, {\Z}(i))$ and $ H^k_{\M} (o_F^S, {\Z}(i))$. For an
integer $n$, the motivic cohomology group with coefficients ${\Z}/n$
is the cohomology of the complex ${\Z}/n(i) = {\Z}(i)\otimes
{\Z}/n$.

For a number field $F$ with ring of integers $o_F$ and a fixed prime
$p$, we let $o^{\prime}_F:=o_F[1/p]$ and simply denote by
$H^k(o^{\prime}_F, {\mathbf{Z}}_p(i))$ the \'{e}tale cohomology
groups
 ${H^{k}_{{\mbox{\'et}}}} (\mathrm{spec} (o^{\prime}_F), {\mathbf{Z}}_p(i))$.
  Soul\'{e}, Dwyer and Friedlander \cite{Soule79, Dwyer-Friedlander85} constructed the  \'{e}tale Chern characters $ch_{i,k}^{(p)}$ which give
   the relation to the algebraic  $K$-theory (tensored with ${\mathbf{Z}}_p$) of $o_F$ and
\'{e}tale cohomology.  They proved that for  all $i\geq 2$ and
$k=1,2$,
$$ch_{i,k}^{(p)} : K_{2i-k}o_F \otimes {\mathbf{Z}}_{p} \rightarrow  H^{k}_{{\mbox{\'et}}} (o^{\prime}_F, {\mathbf{Z}}_p(i))$$
is surjective with finite kernel provided  that $p> 2$. The
Quillen-Lichtenbaum conjecture asserts in fact that $ch_{i,k}^{(p)}$
is an isomorphism for $p> 2$. The odd $p$-primary part  of the
conjecture follows  from the motivic Bloch-Kato Conjecture (cf e.g.
\cite[Theorem 2.7]{Kolster04}) which is
now a theorem thanks to the work of Rost and Voevodsky
\cite{Voevodsky11}. For  $p=2$,  the exact information about the
kernel and the cokernel of the Chern Character $ch_{i,k}^{(2)}$
\cite{Kahn97, Weibel00} has been determined using ideas in
Voevodsky's proof of the Milnor Conjecture \cite{Voevodsky03}. In
particular, $ch_{i,k}^{(2)}$ is an isomorphism   when $2i-k \equiv
0,1,2,7 \pmod 8$ or $F$ is totally imaginary. The link between
$K$-theory and motivic cohomology is given  by motivic Chern
characters \cite{Pushin04}
 $$ch_{i,k}^\M : K_{2i-k} (F)   \longrightarrow  H^k_\M (F, {\Z}(i))$$
which, once tensored with  $\Z_p$, induce $p$-adic Chern characters
\cite[Chapter III]{Levine98}. The Bloch-Kato conjecture implies that this map is
an isomorphism -up to
   $2$-torsion- for all $i\geq 2$ and $k=1,2$.
As a consequence of the Bloch-Kato conjecture, we also have  the
following comparison between motivic cohomology and \'{e}tale
cohomology. If $S$ is a set of primes of $F$ and $o_F ^S$ is the
ring of $S$-integers of $F$,  there are isomorphisms
 \begin{equation}\label{comparison motivic etale}
 H^k_\M (o_F^S, {\Z}(i))  \otimes  {\Z}_{p}  \cong
H^k_\et (o_F^S\left[1/p\right],  {\Z}_{p}(i))
 \end{equation}
for all integers $i\geq 2$, $k=1,2$ and all prime numbers $p$.
%(\cite[section 2]{Kolster04})
For the second cohomology groups, we have  
 $$H^2_\M(o_F^S, {\Z}(i)) \cong \prod_{p} H^2_\et (o_F^S\left[1/p\right],  {\Z}_{p}(i))$$
(cf e.g.\cite{Kolster04}). In fact one can use this property to give
a global model for the \'etale cohomology groups $H^2_\et
(o_F^S\left[1/p\right],  {\Z}_{p}(i))$.  The construction of a
global model in the same way for the  groups $H^1_\et
(o_F^S\left[1/p\right],  {\Z}_{p}(i))$ is more complicated and
performed in  \cite{Chinburgetal98}. Motivic cohomology groups with
finite coefficients  ${\Z}/n(i)$  fit into a long exact sequence
%(\cite [Geisser] hand book of K-theory page 196)
 \begin{equation}\label{motivic long exact sequence}
\cdots  \rightarrow  H^k_\M (X, {\Z}(i)) \stackrel{n}{\rightarrow}
H^k_\M (X, {\Z}(i))
 \rightarrow  H^k_\M (X, {\Z}/n(i))  \rightarrow  H^{k+1}_\M(X, {\Z}(i))
  \rightarrow  \cdots
  \end{equation}
 For a field  $K$, we note that  $H^k_\M (K, {\Z}/n(i))$ coincides with the Galois cohomology group
 $H^k(K, \mu_n^{\otimes i})$.
 One of the  main ingredients in the sequel is the following
   localization  sequence in motivic  cohomology  \cite[Theorem 1.1]{Geisser04}, which relates motivic cohomology
   groups of the ring $o_F^S$ of $S$-integers of $F$ to the motivic cohomology of $F$, and is similar to Soul\'{e}'s
    localization sequence in  \'{e}tale cohomology \cite{Soule79},
 $$ 0 \rightarrow H^2_\M (o_F^S, {\Z}(i))   \rightarrow  H^2_\M (F, {\Z}(i)) \rightarrow
 \oplus_{v\not\in S}  H^1_\M (k_v,  {\Z}(i-1)) \rightarrow 0 $$
  and  isomorphisms  $H^1_\M (o_F^S, {\Z}(i))   \cong  H^1_\M (F, {\Z}(i))$.
Here $k_v$ stands for the residue field at the place $v$. Recall
also that for all  $i \geq 2$
\begin{equation}\label{residual even K-group}
K_{2i-2}k_v \cong H^2_\M (k_v,  {\Z}(i)) = 0
\end{equation}
and
\begin{equation}\label{residual odd K-group}
K_{2i-1}k_v \cong H^1_\M (k_v, \Z(i)) \cong H^0 (k_v,  {\Q/\Z}(i))
\simeq \Z/(q_v^i-1) 
\end{equation}
where  $q_v$  is the order of the residue field  $k_v$.  
Writing the above exact sequence for both  $o_F$  and  $o_F^S$  leads to the exact sequence
 \begin{equation}\label{tame kernel}
0\fd H^2_\M(o_F, {\Z}(i))  \fd   H^2_\M(o_F^S, {\Z}(i)) \fd
 \oplus_{v\in S\setminus S_\infty }  H^1_\M(k_v, {\Z}(i-1) ) \fd 0.
\end{equation}
Using Quillen's localization sequence
%relating the $K$-theory of  $o_F^S$, $F$ and the finite fields $k_v$
$$0 \rightarrow K_{2i-2}  o_F^S \rightarrow  K_{2i-2}F \rightarrow
 \oplus_{v\not\in S}  K_{2i-3}k_v \rightarrow 0 $$
  we obtain the isomorphisms
$$  K_{2i-k}  o_F^S \cong  H^k_\M (o_F^S, {\Z}(i))$$
for  $i\geq 2$, $k=1,2$, up to known  $2$-torsion \cite{ Kahn97,
Weibel00}. Hence in view of  \cite{Tate76, Soule79}, the sequence
(\ref{tame kernel}) is the analogue of the exact sequence
$$0 \rightarrow K_{2}  o_F \rightarrow  K_{2}o_F^S \rightarrow
 \oplus_{v \in S \setminus S_\infty}  k_v^\bullet \rightarrow 0$$
%defining the  tame kernel.
for  $K_2$. 

From the above discussions and the finiteness of  $K_{2i-2}  o_F$
(Borel), it is also clear  that
 % $H^1_\M(o_F^S, {\Z}(i)) \cong  H^1_\M(o_F, {\Z}(i))$ for $i\neq 1$ .
 the motivic  cohomology group  $H^2_\M(o_F, {\Z}(i))$  is finite  for $i\geq 2$.
 We call it the (higher) \emph{motivic tame kernel}.
Its order  $h_i^\M$  is related to the special value of  the
Dedekind $\zeta$-function at the negative integer  $1-i$ (\cite
[$\cdots$]{BurnsGreither03, KNF96}).

Finally fix a finite prime  $v$  in  $F$. Denote by  $\ell$ the
prime integer such that  $v \mid \ell$. Local class-field theory and
Tate duality show  that the $p$-adic Galois cohomology groups $H^k
(F_v, {\mathbf{Z}}_p(i) )$, $k=0,1,2$, are finite for any prime
integer $p \neq \ell$  and trivial for almost all $p$. Define
$$H^k (F_v, {\mathbf{Z}}(i) )^\prime:=  \prod_{p \neq \ell} H^k (F_v, {\mathbf{Z}}_p(i) ).$$
 The groups $H^k (F_v, {\mathbf{Q}}/{\mathbf{Z}}(i)
)^\prime$  and  $H^k (F_v, {\mathbf{Z}}/n(i) )^\prime$ are defined
in the same way. By local duality,
$$
\begin{array}{ccc}
  H^2 (F_v, {\Zp}(i)) & \cong & H^0 (F_v, {\Qp/\Zp}(i-1)) \\
   & \cong & H^0 (k_v, {\Qp/\Zp}(i-1) )
\end{array}
$$
for any prime integer $p \neq \ell$. Therefore
 \begin{equation}\label{p-adic local residual 2/1}
      H^2(F_v, {\Zp}(i) ) \cong H^1(k_v, {\Zp}(i-1) )   \mbox{ once }   v\nmid p.
      \end{equation}
Hence, taking the product over all the primes  $p$  such that
$v\nmid p$:
 \begin{equation}\label{local residual 2/1}
      H^2 (F_v, {\Z}(i) )^\prime \cong H^1_\M(k_v, {\Z}(i-1) ).
      \end{equation}

Furthermore, the local cohomology groups  $H^1 (F_v, {\Zp}(i))$ being
finite for any prime  $p \neq \ell$, we have
$$H^1 (F_v, {\Zp}(i)) = H^0 (F_v, {\Qp/\Zp}(i))  \cong  H^0 (k_v, {\Qp/\Zp}(i))$$
so that
\begin{equation}\label{p-adic local residual 1/1}
      H^1(F_v, {\Zp}(i) ) \cong H^1(k_v, {\Zp}(i) )   \mbox{ once }   v\nmid p
      \end{equation}
and
\begin{equation}\label{local residual 1/1}
      H^1(F_v, {\Z}(i) )^\prime \cong H^1_\M(k_v, {\Z}(i) ).
      \end{equation}

We shall often use the following result which describes Galois
descent and codescent for motivic cohomology groups.

% \cite{Chinburgetal98} based on ideas of Milne and Kahn  (see KO 03 page 334 top of the page).

\begin{thm}  \cite[Theorem 2.13]{Kolster04} \label{galois descent and codescent}
Let  $L/F$ be a finite Galois extension of number fields with group
$G= \Gal(L/F)$. Let $S$ be a finite set of places of $F$ containing the set  $S_\infty$ of infinite places and all those
which ramify  in $L/F$.
For  $i \geq 2$, we have: \\
(i) $H^1_\M(F, {\Z}(i)) \cong H^1_\M(L, {\Z}(i))^G $ and  there are
exact sequences
$$0 \rightarrow H^2_\M(o_L^S, {\Z}(i))_{G} \rightarrow  H^2_\M(o_F^S, {\Z}(i))  \rightarrow   ({\Z}/2)^r \rightarrow 0
\quad  \mbox{for $i$  even}$$
$$0 \rightarrow  ({\Z}/2)^{s_i} \rightarrow H^2_\M(o_L^S, {\Z}(i))_{G} \rightarrow  H^2_\M(o_F^S, {\Z}(i))  \rightarrow   0
\quad  \mbox{for $i$  odd}$$ (ii) For $i$ even,  we have  an exact
sequence
$$ 0 \rightarrow H^1(G, H^1_\M(L, {\Z}(i))) \rightarrow H^2_\M(o_F^S, {\Z}(i))  \rightarrow
H^2_\M(o_L^S, {\Z}(i))^G \rightarrow H^2(G, H^1_\M(L, {\Z}(i)))
\rightarrow 0.$$ (iii) For $i$ odd,   there is an exact sequence
$$ 0 \rightarrow H^1(G, H^1_\M(L, {\Z}(i))) \rightarrow H^2_\M(o_F^S, {\Z}(i))  \rightarrow
H^2_\M(o_L^S, {\Z}(i))^G $$
$$\rightarrow H^2(G, H^1_\M(L, {\Z}(i)))
\rightarrow   ({\Z}/2)^{r-{s_i}} \rightarrow 0.$$ %and
%$$0   \rightarrow   ({\Z}/2)^{s_i}  \rightarrow \hat{H}^{2q-1}(G, H^2_\M(o_L^S, {\Z}(i))) \rightarrow \hat{H}^{2q+1}(G, H^1_\M(L, {\Z}(i)))   \rightarrow 0,$$
%  $$ 0 \rightarrow
%\hat{H}^{2q}(G, H^2_\M(o_L^S, {\Z}(i)))  \rightarrow
%\hat{H}^{2q+2}(G, H^1_\M(L, {\Z}(i)))   \rightarrow
%({\Z}/2)^{r-{s_i}} \rightarrow 0$$ for all $q\in \Z$.
\end{thm}

The above theorem is proven $p$-part by $p$-part. For $p$ odd, the
proof uses in an essential way the fact that the \'etale cohomology
groups $H^k(o_F^S[1/p], \Z_p(i))$ vanish for $k\not=1,2$. When
$p=2$, one replaces \'etale cohomology by  positive \'etale
cohomology \cite[\S 1] {Chinburgetal98} inspired by the ideas from \cite{Kahn93} and \cite{Milne86}.
%Assume that the set $S$ contains the dyadic primes of  $F$.
 Recall the main properties of the positive \'etale cohomology
groups $H^k_+(o_F^S[1/2], \Z_2(i))$. First notice that the positive
\'etale cohomology groups vanish for $k\not=1,2$. In particular, we
have  isomorphisms
 \begin{equation}\label{descent+}
H^1_+(o_F^S[1/2], {\Z}_2(i))  \cong H^1_+(o_L^S[1/2], {\Z}_2(i))^G
\end{equation}
and
 \begin{equation}\label{codescent+}
H^2_+(o_L^S[1/2], {\Z}_2(i))_G \cong  H^2_+(o_F^S[1/2], {\Z}_2(i))
\end{equation}
 as well as the following dimension shifting result:
 \begin{equation}\label{shifting+}
 \hat{H}^{q}(G, H^2_+(o_L^S[1/2], {\Z}_2(i)))  \cong \hat{H}^{q+2}(G, H^1_+(o_L^S[1/2],
{\Z}_2(i)))
\end{equation}
for all $q\in \Z$, given by cup-product \cite[Proposition
3.1]{Chinburgetal98}. We shall use this result to give the
deviation between $ \hat{H}^{q}(G, H^2(o_L^S[1/2], {\Z}_2(i)))$ and
$\hat{H}^{q+2}(G, H^1(o_L^S[1/2],{\Z}_2(i)))$:

%We will use the fact that for a real prime  $v$, the direct sum  \oplus _{ w\mid v} \Z/2(i)$  is  $G_v$-cohomolgica

For  $n \in \Z$ and an infinite prime  $v$, it can readily be seen that 
$$\hat{H}^n(G, \oplus_{w\mid v}\Z_2(i)) = \left\{
\begin{array}{ll} 
 \Z/2 & \mbox{ if }  i+n  \mbox{ is even }, v  \mbox{ is a real prime becoming complex in } L ,\\ 
0 & \mbox{ otherwise } . 
\end{array}
\right.
$$

a)   Suppose first that $i$ is even.
 We have the following short exact sequences \cite[pages 231-232]{Kolster04}
 \begin{equation}\label{H1+H1 i even}
0\rightarrow \oplus _{v\in S_\infty } \oplus_{ w\mid
v}\Z_2(i)\rightarrow H^1_+(o_L^S[1/2], \Z_2(i)) \rightarrow
H^1(o_L^S[1/2], \Z_2(i)) \rightarrow 0
\end{equation}
 and
 \begin{equation}\label{H2+H2 i even}
0\rightarrow  H^2_+(o_L^S[1/2], \Z_2(i)) \rightarrow H^2(o_L^S[1/2],
\Z_2(i)) \rightarrow \oplus _{ v\, real}  \oplus_{ w \, real \mid v}\Z/2(i) \rightarrow 0.
\end{equation}

Since every direct sum  $\oplus_{ w \, real \mid v}\Z/2(i)$  in the exact sequence (\ref{H2+H2 i even})  is $G$-cohomologically trivial (Shapiro), we have 
$$\hat{H}^n(G, H^2_+(o_L^S[1/2], \Z_2(i))) \cong \hat{H}^n(G, H^2(o_L^S[1/2], \Z_2(i)))$$
for  $n \in \Z$. Now, write the cohomology exact sequence corresponding to (\ref{H1+H1 i even}) 

$$  
 0  \rightarrow \hat{H}^{2q+1}(G, H^1_+(o_L^S[1/2], {\Z}_2(i))) 
\rightarrow \hat{H}^{2q+1}(G, H^1(o_L^S[1/2], {\Z}_2(i))) 
\rightarrow ( {\Z}/2)^r
$$
$$
\rightarrow
\hat{H}^{2q+2}(G, H^1_+(o_L^S[1/2], {\Z}_2(i)))  \rightarrow
\hat{H}^{2q+2}(G, H^1(o_L^S[1/2], {\Z}_2(i))) \rightarrow 0.
$$
Therefore by the dimension shifting (\ref{shifting+}) and the above isomorphism, we obtain the exact sequence
 \begin{equation}\label{kolster1}
0 \rightarrow \hat{H}^{2q-1}(G, H^2(o_L^S[1/2], {\Z}_2(i)))
\rightarrow \hat{H}^{2q+1}(G, H^1(o_L^S[1/2], {\Z}_2(i)))
\rightarrow ( {\Z}/2)^r
$$
$$\rightarrow
\hat{H}^{2q}(G, H^2(o_L^S[1/2], {\Z}_2(i)))  \rightarrow
\hat{H}^{2q+2}(G, H^1(o_L^S[1/2], {\Z}_2(i))) \rightarrow 0
\end{equation}
for all $q\in \Z$. \\

b) Assume now that   $i$ is odd.  Then we have a six-term exact sequence
\cite[page 231]{Kolster04}
 \begin{equation}\label{kolster2}
0\rightarrow \oplus _{ w\, complex} \Z_2(i)\rightarrow
H^1_+(o_L^S[1/2], \Z_2(i)) \rightarrow H^1(o_L^S[1/2], \Z_2(i))
\rightarrow  \oplus _{ w\, real} {\Z}/2$$
$$
\rightarrow H^2_+(o_L^S[1/2], \Z_2(i)) \rightarrow H^2(o_L^S[1/2],
\Z_2(i)) \rightarrow 0
\end{equation}
In this exact sequence  the map
$$H^1(o_L^S[1/2], {\Z_2}(i))   \rightarrow    \oplus _{ w\, real} {\Z}/2 = ({\Z}/2)^{r_1(L)}$$
factors through $H^1(o_L^S[1/2], {\Z_2}(i))/2$ and induces the
signature map  $\mathrm{sgn}_L$. 
Split the above exact sequence into four short ones : 
$$ 0\rightarrow \oplus _{ w\, complex} \Z_2(i)\rightarrow H^1_+(o_L^S[1/2], \Z_2(i)) \rightarrow C \rightarrow 0 $$
$$0\rightarrow C  \rightarrow H^1(o_L^S[1/2], \Z_2(i)) \rightarrow  ({\Z}/2)^{r_1(L)-\delta_i(L)}\rightarrow  0 $$
$$0 \rightarrow  ({\Z}/2)^{r_1(L)-\delta_i(L)}\rightarrow  ({\Z}/2)^{r_1(L)} \rightarrow  ({\Z}/2)^{\delta_i(L)} \rightarrow 0 $$
\begin{equation}
\label{s_i as kernel}
0 \rightarrow  ({\Z}/2)^{\delta_i(L)} \rightarrow  H^2_+(o_L^S[1/2], \Z_2(i)) \rightarrow H^2(o_L^S[1/2], \Z_2(i)) \rightarrow 0
\end{equation}
where obviously  $C$  contains   $2H^1(o_L^S[1/2], \Z_2(i))$  and  the factor group  $C/2H^1(o_L^S[1/2], \Z_2(i))$  is the kernel of the signature map  
$$\mathrm{sgn}_L : H^1(L, {\Z_2}(i))/2 \longrightarrow  ({\Z}/2)^{r_1(L)}$$ 
introduced in section \ref{subsection Signature}. 

Assume   
$$\hat{H}^n(G, ({\Z}/2)^{\delta_i(L)})=0   \mbox{ for all } n \in \Z.$$ 
Then  $H^2_+(o_L^S[1/2], \Z_2(i))$  and  $H^2(o_L^S[1/2], \Z_2(i))$  have the same  $G$-cohomology.  
Also, by Shapiro's Lemma  
$\hat{H}^n(G, ({\Z}/2)^{r_1(L)})=0   \mbox{ for all } n \in \Z$.   
Therefore 
$\hat{H}^n(G, ({\Z}/2)^{r_1(L)-\delta_i(L)})=0   \mbox{ for all } n \in \Z$   
and the two  $G$-modules  $C$  and  $H^1(o_L^S[1/2], \Z_2(i))$  also have the same  $G$-cohomology.  
Taking the  cohomology of the first short exact sequence above and using the dimension shifting (\ref{shifting+}) yields 
the following exact sequence   
 \begin{equation}\label{kolster3}
0 \rightarrow \hat{H}^{2q-2}(G, H^2(o_L^S[1/2], {\Z}_2(i)))
\rightarrow \hat{H}^{2q}(G, H^1(o_L^S[1/2], {\Z}_2(i))) \rightarrow
( {\Z}/2)^r $$
$$\rightarrow
\hat{H}^{2q-1}(G, H^2(o_L^S[1/2], {\Z}_2(i)))  \rightarrow
\hat{H}^{2q+1}(G, H^1(o_L^S[1/2], {\Z}_2(i))) \rightarrow 0
\end{equation}
for all $q\in \Z$ under the hypothesis $\hat{H}^n(G, ({\Z}/2)^{\delta_i(L)})=0   \mbox{ for all } n \in \Z.$ 

 When $q=0$ (and  $i$ odd), we have the above exact sequence  (\ref{kolster3}) as soon as hypothesis  ($ \mathcal{H}_i $)  is fulfilled. Furthermore, we have a  precise description of the image of the
 map
 $$\hat{H}^{0}(G, H^1(o_L^S[1/2], {\Z}_2(i)))   \rightarrow  ( {\Z}/2)^r.$$
 Indeed, the signature map  gives rise to the
following commutative diagram
$$
  \begin{array}{cccccc}
     H^1(o_F^S[1/2], {\Z}_2(i))   &    &   \rightarrow   &  ( {\Z}/2)^r  &
     \\

    {\downarrow}                    &                       &                  &  {\downarrow}\wr &
      \\
       \hat{H}^0(G, H^1(o_L^S[1/2], {\Z}_2(i)))   &     &   \rightarrow   &  ( {\Z}/2)^r  &
     \\
  \end{array}
  $$
where the  left hand vertical map is surjective (Theorem \ref{galois descent and codescent} (i)).
An easy diagram chasing  yields 
$$\coker (\hat{H}^0(G, H^1(o_L^S[1/2], {\Z}_2(i)))     \rightarrow (
{\Z}/2)^r) \cong  ( {\Z}/2)^{s_i}$$
so that the beginning of the exact sequence  (\ref{kolster3}) for  $q=0$ (and  $i$ odd)  is written  
 \begin{equation}\label{kolster4}
0 \rightarrow H_1(G, H^2(o_L^S[1/2], {\Z}_2(i)))
\rightarrow \hat{H}^0(G, H^1(o_L^S[1/2], {\Z}_2(i))) \rightarrow ({\Z}/2)^r   \rightarrow ({\Z}/2)^{s_i} \rightarrow 0  
\end{equation}
under  hypothesis ($ \mathcal{H}_i $).  
%\begin{rem}\label{usual sequence}
%The only assumption we need to prove (\ref{kolster4} ) is the vanishing of  $\hat{H}^{-1}(G, ({\Z}/2)^{\delta_i(L)})$.
%In fact, the hypothesis ($ \mathcal{H}_i $) in the sequel can be replaced by the nullity of  $\hat{H}^n(G, ({\Z}/2)^{\delta_i(L)})$.
%\end{rem}

\begin{rem}\label{kcyclic shifting}
If the extension $L/F$ is unramified at infinite prime, Theorem  \ref{galois descent and codescent} shows that $H^2_\M( L, {\Z}(i) )$ satisfy Galois co-descent. Il follows that the kernel and  cokernel of the functorial map 
$$f_i : H^2_\M( o_F^S, {\Z}(i) ) \rightarrow
H^2_\M( o_L^S, {\Z}(i) )^G$$
 are given by  $\ker (f_i) \cong \hat{H}^{-1} (G, H^2_\M( o_L^S {\Z}(i) )  \cong {H}^{1} (G, H^1_\M( o_L^S, {\Z}(i) )$ and  $\coker (f_i) \cong \hat{H}^{0} (G, H^2_\M( o_L^S, {\Z}(i)   \cong {H}^{2} (G, H^1_\M( o_L^S {\Z}(i) ))$.
If, Moreover, $G$ is a cyclic group, we then have the dimension shifting : 
\begin{equation}
 \hat{H}^{q}(G, H^2_\M(o_L^S[1/2], {\Z}(i)))  \cong \hat{H}^{q+2}(G, H^1_\M(o_L^S[1/2],
{\Z}(i)))
\end{equation}
for all $i\geq 2$ and all $q\in \Z$, given by cup-product.
\end{rem}

Notice that similar results hold in the local case.  We mention the
 following isomorphisms  which will be used in the next
sections:

For any local extension $L_w/F_v$ with Galois group $G_v$,  denote
by $k_v$ (resp. $k_w$) the residue field of $F_v$ (resp.$L_w$). The
group  $G_v$  acts on  $H^1_\M(k_w, {\Z}(i))$  via the natural
composite map   $G_v  \twoheadrightarrow G_v/I_v \stackrel{\sim}
\rightarrow \Gal(k_w/k_v)$.
 It follows from the isomorphisms (\ref{local residual 2/1}) and (\ref{local residual 1/1})  that
  %({\bf voir d'autres references plut�t que passer par les corps locaux ....})
\begin{equation}\label{residual descent}
 H^1_\M(k_v, {\Z}(i)) \cong H^1_\M(k_w, {\Z}(i))^{G_v},
  \end{equation}
 \begin{equation}\label{residual codescent}
 H^1_\M(k_w, {\Z}(i))_{G_v} \cong H^1_\M(k_v, {\Z}(i))
  \end{equation}
 and
 \begin{equation}\label{residual shifting}
 \hat{H}^q(G_v,  H^1_\M(k_w, {\Z}(i-1))) \cong \hat{H}^{q+2}(G_v, H^1_\M(k_w, {\Z}(i))),
  \end{equation}
 for all $q\in  {\Z}$ and $i\geq 2$.
%The exact sequence (\ref{motivic long exact sequence}) then shows that the natural "transfer" map
%$$\mathrm{tr}_i : H^2( o_{L}, {\Z}(i) )_{G} \rightarrow
%H^2( o_{F}, {\Z}(i) $$
%is surjective, up to known $2$-torsion.  In the next section we work out the kernel of this map explicitly for any Galois
%extension  $L/F$.

%\label{sec:greetings} Hello! \section{Referencing} I greeted in section~\ref{sec:greetings}.

\section {Genus formula for motivic tame kernels}\label{Genus formula for motivic tame kernels}
Let $L/F$ be a finite Galois extension of algebraic number fields with Galois group $G$. 
From here on, unless specified otherwise, the set  $S$  consists of infinite primes and those which ramify in  $L/F$. 
The short exact sequence  (\ref{tame kernel}) yields the following commutative diagram:
$$
\begin{array}{cccccccccccc}
   &  &    &
&H^2_\M(o_L, {\Z}(i))_G & {\rightarrow} & H^2_\M(o_L^S, {\Z}(i))_G &
\rightarrow &
({\oplus_{v \in {S \setminus S_\infty}}}({\oplus_{w \mid v}} {H^1_\M(k_w,\Z(i-1))))_G} & \fd & 0\\
           &  & & &  \mathrm{tr}_i \downarrow &    &   \downarrow &             & \wr   \downarrow \\
     & &  0 & \rightarrow & H^2_\M(o_F, {\Z}(i)) &
{\rightarrow} & H^2_\M(o_F^S, {\Z}(i))& \rightarrow & \oplus_{v \in
{S\setminus S_\infty}} {H^1_\M(k_v,\Z(i-1))} & \rightarrow &  0
  \end{array}
  $$
 where the right hand vertical map is an isomorphism by (\ref{residual codescent}). Let
$${\alpha} :  H_1(G, H^2_\M(o_L^S, {\Z}(i))) {\rightarrow} {H}_1(G, \oplus_{v\in S\setminus S_\infty, w\mid v}H^1_\M(k_w,\Z(i)))$$
be  the map given by the homology of the exact sequence (\ref{tame
kernel}). Using Theorem \ref{galois descent and codescent} ( i), we then derive the following
\begin{proposition}\label{beginning genus}
If  $i\geq 2$ is even, we have an exact sequence
$$0 \rightarrow  {\mathrm{coker}} (\alpha)
\rightarrow H^2_\M(o_L, {\Z}(i))_{G} {\rightarrow} H^2_\M(o_F,
{\Z}(i)) {\rightarrow} ({\Z}/2)^r \rightarrow 0.$$ If $i\geq 3$ is
odd, the map $ \mathrm{tr}_i :  H^2_\M( o_{L}, {\Z}(i) )_{G}
\rightarrow H^2_\M( o_{F}, {\Z}(i)) $ is surjective and its kernel
fits into an exact sequence
$$0 \rightarrow  {\mathrm{coker}} \alpha
\rightarrow \ker  \mathrm{tr}_i {\rightarrow} ( {\Z}/2)^{s_i}
\rightarrow 0.$$
  In particular
$$  \frac{|H^2_\M(o_L, {\Z}(i))_{G}|}{|H^2_\M(o_F, {\Z}(i))|}
 =2 ^{n_i} {|\mathrm{coker}(\alpha)|}$$
 %= 2 ^{n_i} \frac{\prod_{v\in S\setminus S_\infty} |H_1(G_v, H^1_\M(k_w, {\Z}(i-1) ))  |}{|\mathrm{Im}(\alpha)|}.$$
 where $n_i =-r$ if $i$ is even and $n_i =s_i$ if $i$ is odd.
\qed
\end{proposition} 

For each prime  $v \in {S\setminus S_\infty}$, we fix a prime  $w$
of  $L$ above  $v$  and denote by  $G_v$  the local Galois group
$\Gal(L_w/F_v)$. We first express the  homology groups
$H_1(G_v,H^1_\M(k_w, {\Z}(i-1) )) $ in terms of the ramification in
$L/F$. Denote by  $q_v$  the order of the residue field  $k_v$  and
by  $e_v$ the ramification index of $v$ in the extension $L/F$. If
$\ell$ is the residue characteristic of $F_v$, we denote by  $
e_v^\prime$ the order of the non-$\ell$-primary part of the inertia
group of $v$ in the extension $L/F$.

\begin{lem}\label{lemma evi}
For each  finite prime  $v$, we have
$$H_1(G_v,H^1_\M(k_w, {\Z}(i-1) ) ) \simeq \Z/{e_{v}^{(i)}}$$
where  $e_v^{(i)} :=gcd (e_v, q_v^i-1)$. In particular, when  $G_v$
is abelian, then  $e_v^{(i)}=e_v^\prime$ and therefore does not
depend on the integer $i$.
\end{lem}

{\bf Proof} By the dimension shifting  (\ref{residual shifting}), we
have
$$ {H}_1(G_v,  H^1_\M(k_w, {\Z}(i-1))) \cong \hat{H}^{0}(G_v, H^1_\M(k_w, {\Z}(i)))
  .$$
 Let $E_w=F_{v}^{nr}\cap L_w$ be the maximal unramified extension of $F_v$ contained in $L_w$ and  $I_v:= \Gal(L_w/E_w)$ be the inertia group.

Since $I_v$  acts trivially on  $H^1_\M(k_w,\Z(i))$, we have
$$
  \begin{array}{lll}
      \hat{H}^0(G_v, H^1_\M(k_w,\Z(i))) & \cong &  H^1_\M(k_v,\Z(i))/ N_{G_v/I_v}(N_{I_v} H^1_\M(k_w,\Z(i))) \\
        & \cong & H^1_\M(k_v,\Z(i))/e_v N_{G_v/I_v}( H^1_\M(k_w,\Z(i)))
     \end{array}
$$
where  $N_{G_v/I_v}$ (resp.  $N_{I_v}$) is the norm map
corresponding to the extension  $E_w/F_v$ (resp.  $L_w/E_w$). By
\cite[section 12, page 185 remark]{Quillen72}, the norm map
$$N_{G_v/I_v}: H^1_\M(k_w,\Z(i))  \rightarrow   H^1_\M(k_v,\Z(i))$$
is surjective. Hence
$$
  \begin{array}{lll}
\hat{H}^0(G_v, H^1_\M(k_w,\Z(i))) & \simeq &  H^1_\M(k_v,\Z(i))/e_v\\
       &  \simeq &  \Z/{e_{v}^{(i)}}
     \end{array}
$$
  as required.
The last assertion comes from the fact that  $I_v/I_{v,1}$ is
isomorphic to a quotient of  $k_v^\bullet$, where  $I_{v,1}$  is the
first ramification group.
  \qed

For each prime number  $p$, and each finite non-$p$-adic prime  $v$
of  $F$, we consider the composite map
$$
\begin{array}{cccc}
  H^1_\M(F, {\Z}(i)) \otimes \Z_p &   \stackrel{\sim} \longrightarrow  &  H^1 (F, {\Z_p}(i)) &  (\mbox{by the isomorphism (\ref{comparison motivic etale}}))   \\
  &  \stackrel{res_v} \longrightarrow &   H^1 (F_v, {\Z_p}(i)) & \\
  & \stackrel{\sim} \longrightarrow  &   H^1 (k_v, {\Z_p}(i)) &   (\mbox{by the isomorphism (\ref{p-adic local residual 1/1}}))    \\
 & \stackrel{\sim} \longrightarrow  &   H^1_\M (k_v, {\Z}(i)) \otimes \Z_p & (\mbox{by the isomorphism (\ref{residual odd K-group}})).  \\

\end{array}
$$
%(voir si on ne peut pas mettre dans un meme cadre le premier et le
%dernier isomorphisme)

Let  $p^m$  be the exact power of  $p$  dividing the order $q_v^i-1$
of $H^1_\M(k_v, {\Z}(i))$. We then have a canonical map
\begin{equation*}
H^1_\M(F, {\Z}(i))/p^m{\rightarrow} H^1_\M (k_v, {\Z}(i))/p^m.
\end{equation*}
Taking the direct sum over all primes $p\mid q_v^i-1$  yields a map
$$  H^1_\M(F, {\Z}(i))/(q_v^i-1){\rightarrow}  H^1_\M(k_v, {\Z}(i)).$$
Finally, composing with the canonical surjection
$$  H^1_\M(F, {\Z}(i)) {\rightarrow}  H^1_\M(F, {\Z}(i))/(q_v^i-1)$$
we have a map 
\begin{equation}\label{phiv}
\varphi_v : H^1_\M(F, {\Z}(i)) {\rightarrow}  H^1_\M(k_v, {\Z}(i)).
\end{equation}
Now, introduce the map 
\begin{equation}\label{phiS}
\beta_S : \hat{H}^0(G, H^1_\M(L,\Z(i))) \rightarrow  
\hat{H}^0(G, \oplus_{w \mid v \in S\setminus S_\infty}H^1_\M(k_w,\Z(i)))
 \end{equation}
given by taking the cohomology on the direct sum
$$\oplus_{w \mid v \in S\setminus S_\infty} \varphi_w :
H^1_\M(L, {\Z}(i)) {\rightarrow}  \oplus_{w \mid v \in S\setminus
S_\infty}H^1_\M(k_w, {\Z}(i)).$$

The exact sequences (\ref{kolster1}) ($i$ even), (\ref{kolster4})
($i$ odd, under hypothesis  ($\mathcal{H}_i$)) and  the shifting (\ref{residual shifting}) give rise to
the commutative diagram
 \begin{equation}\label{diagram alpha betaS}
 \begin{array}{ccc}
     {H}_1(G, H^2_\M(o_L^S,\Z(i)))& \stackrel{\alpha} \longrightarrow & {H}_1(G, \oplus_{w \mid v \in S\setminus S_\infty}H^1_\M(k_w,\Z(i-1))) \\
   \downarrow &  & \downarrow \wr\\
    \hat{H}^0(G, H^1_\M(L,\Z(i))) & \stackrel{\beta_S} \longrightarrow & \hat{H}^0(G, \oplus_{w \mid v \in S\setminus S_\infty}H^1_\M(k_w,\Z(i)))
 \end{array}
  \end{equation}
Hence: \\
 (i) for $i$  even, the  exact sequence  (\ref{kolster1}) implies that the left vertical map is surjective so that
 $$\mathrm{coker}(\alpha) \cong \mathrm{coker}(\beta_S)$$
 (ii) for $i$ odd, using  (\ref{kolster4}), we have the following exact sequence
 $$0 \rightarrow \mathrm{ker}(\alpha) \rightarrow \mathrm{ker}(\beta_S) \rightarrow (\Z/2)^{r-s_i} \rightarrow \mathrm{coker}(\alpha) \rightarrow \mathrm{coker}(\beta_S) \rightarrow 0$$ 
 under hypothesis  ($\mathcal{H}_i$).  
%Let  $n$ be the exponent  of the Galois group $G=\Gal(L/F)$.
Since $H^1_\M(L,\Z(i))$ verifies Galois descent (Theorem \ref{galois descent and codescent}(i)),
the above map  $\beta_S$ can be written, by Shapiro's lemma, as
$$ H^1_\M(F,\Z(i))/N_G(H^1_\M(L,\Z(i))) \longrightarrow   \oplus_{v\in S\setminus S_\infty} H^1_\M(k_v,\Z(i))/N_{G_v}(H^1_\M(k_w,\Z(i))).$$

\begin{defi}\label{def N} 
We define the normic subgroup  $H^{1, \mathcal{N}}_\M(F, {\Z}(i))$  to be the set of elements 
$a \in H^1_\M(F, {\Z}(i))$ whose image by  $\varphi_v$  belongs to
 $N_{L_w/F_v}(H^1_\M(k_w, {\Z}(i)))$,  for all $w\mid v \in S \setminus S_\infty$.
In other words, $H^{1, \mathcal{N}}_\M(F, {\Z}(i))$ is the kernel of
the map
$$H^1_\M(F,\Z(i)) \longrightarrow
\oplus_{v\in S\setminus S_\infty}
H^1_\M(k_v,\Z(i))/N_{G_v}(H^1_\M(k_w,\Z(i)))$$ induced by the maps
$\varphi_v$.
\end{defi}
 Hence
$$\mathrm{Im}(\beta_S) \cong H^1_\M(F, {\Z}(i))/ H^{1, \mathcal{N}}_\M(F,
{\Z}(i)).$$

Now by Lemma 2.2, it follows that:
$$\mid \mathrm{coker}(\beta_S) \mid = \frac{\prod_{v\in S\setminus S_\infty}e_{v}^{(i)}}
{[H^1_\M(F, {\Z}(i)): H^{1, \mathcal{N}}_\M(F, {\Z}(i))]}. $$ Taking
into account the above discussion, Proposition \ref{beginning genus}
yields the following genus formula:
\begin{thm}\label{motivic genus}
Let  $L/F$ be a finite Galois extension of number fields with Galois
group  $G$  and  $i\geq 2$. When  $i$  is odd, we assume hypothesis  ($\mathcal{H}_i$) holds.  
Denote by  $S$  the set of infinite
primes and those which ramify in  $L/F$. Then
$$\frac{|H^2_\M(o_L, {\Z}(i))_{G}|}{|H^2_\M(o_F, {\Z}(i))|}= 2 ^{\nu_i}
\frac{\prod_{v \in S \setminus S_\infty}e_{v}^{(i)}}{[H^1_\M(F,
{\Z}(i)): H^{1, \mathcal{N}}_\M(F, {\Z}(i))]}$$ where for each
finite prime  $v$, $e_{v}^{(i)} :=gcd (e_v, q_v^i-1)$,  $e_v$ being
the ramification index of  $v$   in the extension  $L/F$,  $\nu_i
=-r$ for $i$ even and $s_i\leq \nu_i \leq r$ for $i$ odd.
\qed 
\end{thm}

\begin{rem}
 When  $L/F$ is unramified at infinite primes, the power $2^{\nu_i}$ disappears.   In  the general case, the quantity $\nu_i $, $i$ odd,
  may take any of the two bounds as the following examples show.
  \begin{enumerate}
\item  Suppose that $L/F$ is unramified at all finite primes. Then  $ \mathrm{coker}(\alpha) \cong \mathrm{coker}(\beta_S)=0$ hence $\nu_i =s_i$.

\item  Let $L/F$ be a CM-extension such that the natural map
$$H^1_\M(F,\Z(i))/2  \longrightarrow   \oplus_{v\in S\setminus S_\infty} H^1_\M(k_v,\Z(i))/2$$
is an isomorphism. Then
$$\beta_S : \hat{H}^0(G, H^1_\M(L,\Z(i)))  \rightarrow \hat{H}^0(G,  \oplus_{v\in S\setminus S_\infty, w\mid v}H^1_\M(k_w,\Z(i)))$$
is surjective and is in fact an isomorphism taking into account the
orders. An easy diagram chasing in (\ref{diagram alpha betaS}) shows
that, in this situation:

$$\mathrm{coker}(\alpha)  \simeq  (\Z/2)^{r-s_i} $$
 so that   $\nu_i =r  >  0$ (Proposition \ref{beginning genus}).
As  an  explicit example consider the case $F = \Q$. The nullity of
the \'etale cohomology group  ($i$  is odd)
    $H^2_\et (o^\prime_F,  {\Z}_{2}(i))$
and Kummer theory show that
$$H^1_\et (F,  {\Z}_{2}(i))/2 \cong H^1_\et (o^\prime_F,  {\Z}/2(i))
\cong  U^\prime_F/U^{\prime 2}_F = <-{\overline 1},{\overline 2}> $$
where   $o^\prime_F = \Z[1/2]$  and   $U^\prime_F$  denotes the
group of units in  $o^\prime_F$. In this case, the map
$$ U^\prime_F/U^{\prime 2}_F   \rightarrow  k_3^\bullet/k_3^{\bullet 2} \oplus k_5^\bullet/k_5^{\bullet 2}$$
is an isomorphism and all we have to do is to take  $L =
\Q(\sqrt{-15})$. Here, the signature map
$$H^1_\et (F,  {\Z}_{2}(i))/2  \rightarrow \Z/2$$
is surjective and $s_i=0$. Hence  $\nu_i =r \neq s_i$.
   \end{enumerate}
\end{rem}

To remove the ambiguity on the $2$-power in the above theorem when
$i$ is odd, we use positive cohomology. Start with the short exact
sequence
$$0\fd H^2_+ (o_F[1/2], {\Z_2}(i))  \fd   H^2_+ (o_F^S[1/2], {\Z_2}(i)) \fd  \oplus_{v}  H^2(F_v, {\Z_2}(i) ) \fd 0$$
where  $v$  runs over finite non-dyadic primes in  $S$ \cite[page
336]{Kolster03}. According to the isomorphism (\ref {p-adic local
residual 2/1})  and the vanishing of  $H^1(k_v, {\Z_2}(i-1))$  for a
dyadic prime  $v$ (the isomorphism (\ref{residual odd K-group})
tensored with  $\Z_2$), we obtain the exact sequence
 $$0\fd H^2_+ (o_F[1/2], {\Z_2}(i))  \fd   H^2_+ (o_F^S[1/2], {\Z_2}(i)) \fd
 \oplus_{v\in S\setminus S_\infty }  H^1(k_v, {\Z_2}(i-1) ) \fd 0. $$
Hence the following commutative diagram:
$$
\begin{array}{cccccccccccc}
       &
&H^2_+ (o_L[1/2], {\Z_2}(i))_G & {\rightarrow} & H^2_+ (o_L^S[1/2],
{\Z_2}(i))_G & \rightarrow &
({\oplus_{v \in {S \setminus S_\infty,w \mid v}}} {H^1(k_w,\Z_2(i-1))))_G} & \fd & 0\\
            & &  \mathrm{tr}_i ^+\downarrow &    &  \wr  \downarrow &             & \wr   \downarrow \\
        & & H^2_+ (o_F[1/2], {\Z_2}(i)) &
{\hookrightarrow} & H^2_+ (o_F^S[1/2], {\Z_2}(i))& \rightarrow &
\oplus_{v \in {S\setminus S_\infty}} {H^1(k_v,\Z_2(i-1))} &
\rightarrow &  0
  \end{array}
  $$
 which shows that the left vertical map $\mathrm{tr}_i ^+$ is surjective and that
$$ \ker (\mathrm{tr}_i ^+) \simeq {{\mathrm{coker}}
(\alpha}^+ :  H_1(G, H^2_+ (o_L^S[1/2], {\Z_2}(i))) {\rightarrow}
{H}_1(G, \oplus_{w \mid v \in S\setminus
S_\infty}H^1(k_w,\Z_2(i-1)))).$$ Now the dimension shifting
(\ref{shifting+}) shows that
$$\mathrm{coker} (\alpha^+) \simeq \mathrm{coker} (\beta_S ^+ :  \hat{H}^0(G, H^1_+ (L, {\Z_2}(i)))
\rightarrow  \hat{H}^0(G, \oplus_{w\mid v\in S\setminus
S_\infty}H^1(k_w,\Z_2(i)))$$ where  $\beta_S ^+$  is defined in the
same way as  $\beta_S$  with positive cohomology instead of motivic
cohomology. More precisely, for each finite prime $v$,  let
$$\varphi_v^+ : H^1_+ (F, {\Z_2}(i))) {\rightarrow}  H^1(k_v,\Z_2(i)) $$
be the map obtained by composing the map
$$H^1_+(F, \Z_2(i)) \rightarrow H^1(F, \Z_2(i))$$
(in the exact sequences (\ref{H1+H1 i even})  and  (\ref{kolster2}))
with  the map  $\varphi_v$  introduced in  (\ref{phiv}) tensored
with  $ \Z_2$. The map  $\beta_S ^+$  is defined by considering the
direct sum
$$\oplus_{w\mid v\in S\setminus S_\infty} \varphi_w^+ : H^1_+ (L, {\Z_2}(i))) {\rightarrow}
\oplus_{w\mid v\in S\setminus S_\infty} H^1(k_w,\Z_2(i)) $$ 
and passing to the $G$-cohomology. 
In this way, we get the following  $+$-analog of definition  \ref{def N}. 

\begin{defi}\label{def N+} 
We define the (plus) normic subgroup $H^{1, \mathcal{N}}_+(F,{\Z_2}(i))$  to be the kernel of the map
$$H^1_+(F,\Z_2(i)) \longrightarrow \oplus_{v\in S\setminus S_\infty}
H^1(k_v,\Z_2(i))/N_{G_v}(H^1(k_w,\Z_2(i)))$$ induced by the maps  $\varphi_v^+$.
\end{defi}

Since $H^1_+ (o_L^S[1/2], {\Z_2}(i))$ verifies Galois descent by
(\ref{descent+}), it follows exactly as in the proof of  Theorem
\ref{motivic genus} that
$$\mid \mathrm{coker}(\beta_S^+) \mid = \frac{\prod_{v\in S\setminus S_\infty}2^{n_v^{(i)}}}{[H^1_+(F, {\Z_2}(i)): H^{1,\mathcal{N}}_+(F, {\Z_2}(i))]}$$
where for $v\in S\setminus S_\infty$,  $2^{n_v^{(i)}}$ is the exact
power of $2$ dividing $e_{v}^{(i)}$. We summarize
\begin{thm}\label{genus+}
Let  $L/F$ be a finite Galois extension of number fields with Galois
group  $G$  and  $i\geq 2$. Denote by  $S$  the set of infinite
primes and those which ramify in  $L/F$. Then the natural map
$$\mathrm{tr}_i ^+ : H^2_+(o_L[1/2], {\Z_2}(i))_{G} \rightarrow H^2_+(o_F[1/2], {\Z_2}(i))$$
is surjective and we have the following genus formula
$$\frac{|H^2_+(o_L[1/2], {\Z_2}(i))_{G}|}{|H^2_+(o_F[1/2], {\Z_2}(i))|}=
\frac{\prod_{v\in S\setminus S_\infty}2^{n_v^{(i)}}}{[H^1_+(F,
{\Z_2}(i)): H^{1,\mathcal{N}}_+(F, {\Z_2}(i))]} $$ where for $v\in
S\setminus S_\infty$,  $2^{n_v^{(i)}}$ is the exact power of $2$
dividing $e_{v}^{(i)}$. 
\qed 
\end{thm}

We also need the following comparison result:

\begin{proposition}\label{compare + motivic}
Let  $L/F$ be a finite Galois extension of number fields with Galois
group  $G$  and  $i\geq 2$. When  $i$  is odd, we assume hypothesis  ($\mathcal{H}_i$) holds. Then
$$\frac{|H^2_+(o_L[1/2], {\Z_2}(i))_{G}|}{|H^2_+(o_F[1/2], {\Z_2}(i))|}= 2^{u_i} \frac{|H^2(o_L[1/2], {\Z_2}(i))_{G}|}{|H^2(o_F[1/2], {\Z_2}(i))|}$$
where  $u_i=r$  is the number of the real primes of  $F$  ramified in $L$  for $i$ even\\
and  $u_i=-s_i =- s_i(L/F)$  for $i$ odd.
\end{proposition}

{\bf Proof} Denote by  $S$  the set of infinite primes and those
which ramify in  $L/F$. Assume first that  $i$  is odd. In this case
we have a commutative diagram (cf the exact sequences
(\ref{kolster2}) and  (\ref{kolster4}))
$$
\begin{array}{cccccccccccc}
   &  &    &
&((\Z/2)^{\delta_i(L)})_G & {\rightarrow} & H^2_+ (o_L[1/2], {\Z_2}(i))_G
& \rightarrow &
H^2(o_L[1/2], {\Z_2}(i))_G  & \fd & 0\\
           &  & & &   \downarrow &    &   \mathrm{tr}^+_i \downarrow &             & \mathrm{tr}_i   \downarrow \\
     & &  0 & \rightarrow & (\Z/2)^{\delta_i(F)} &
{\rightarrow} & H^2_+ (o_F[1/2], {\Z_2}(i))& \rightarrow & H^2(o_L[1/2], {\Z_2}(i)) & \rightarrow &  0\\
          &  & & &   \downarrow &    &   &             &     \\
          &  & & &   (\Z/2)^{s_i(L/F)} &    &    &             &  \\
          &  & & &   \downarrow &    &   &             &     \\

          &  & & &   0 &    &   &             &
  \end{array}
  $$

with exact lines and columns. 
The left hand vertical map is injective since its kernel  $\hat{H}^{-1}(G, (\Z/2)^{\delta_i(L)})$  vanishes by hypothesis  ($\mathcal{H}_i$).   Since the maps $\mathrm{tr}^+_i
$, $\mathrm{tr}_i$ are surjective according to Theorem \ref{genus+}
and  Proposition \ref{beginning genus} respectively, we have a short
exact sequence 
\begin{equation}\label{tr tr+}
0 \rightarrow \ker (\mathrm{tr}_i ^+) \rightarrow \ker (\mathrm{tr}_i) \rightarrow (\Z/2)^{s_i(L/F)} \rightarrow 0.  
\end{equation}
In particular
$$\mid \ker (\mathrm{tr}_i) \mid  = 2^{s_i}  \mid \ker (\mathrm{tr}_i ^+)\mid .$$
If instead  $i$  is even, the exact sequence  (\ref{H2+H2 i even})
yields the following commutative diagram
$$
\begin{array}{cccccccccccc}
   & & 0  & \rightarrow & H^2_+ (o_L[1/2], {\Z_2}(i)) _G &
{\rightarrow} & H^2(o_L[1/2], {\Z_2}(i)) _G& \rightarrow & (\oplus_{w \, real} \Z/2)_G & \rightarrow &  0\\
           &  & & & \mathrm{tr}^+_i \downarrow   &    &   \mathrm{tr}_i   \downarrow &             & \downarrow\\
     & &  0 & \rightarrow & H^2_+ (o_F[1/2], {\Z_2}(i))  &
{\rightarrow} & H^2(o_F[1/2], {\Z_2}(i)) & \rightarrow & \oplus_{v \, real} \Z/2 & \rightarrow &  0\\
&  & & &   &    &    &             & \downarrow\\
&  & & &   &    &    &             & (\Z/2)^{r(L/F)}\\
&  & & &   &    &    &             & \downarrow\\
&  & & &   &    &    &             & 0

  \end{array}
$$
By Shapiro's lemma,  $H_1(G, \oplus_{w \, real} \Z/2)=0$, hence the
exactness on the left hand side of the top sequence. The right
vertical map being injective, we have
$$ \ker (\mathrm{tr}_i) \cong \ker (\mathrm{tr}_i ^+)$$
and
$$ \coker (\mathrm{tr}_i) \cong (\Z/2)^{r(L/F)}. \qed $$
Combining the above proposition with Theorem \ref{genus+} leads to

\begin{cor}\label{odd 2-primary genus}
Let  $L/F$ be a finite Galois extension of number fields with Galois
group  $G$  and  $i\geq 2$ odd. Assume hypothesis  ($\mathcal{H}_i$) holds. Denote by  $S$  the set of infinite
primes and those which ramify in  $L/F$. Then
$$\frac{|H^2(o_L[1/2], {\Z_2}(i))_{G}|}{|H^2(o_F[1/2], {\Z_2}(i))|}= 2 ^{s_i}
\frac{\prod_{v\in S\setminus S_\infty}2^{n_v^{(i)}}}{[H^1_+(F,
{\Z_2}(i)): H^{1,\mathcal{N}}_+(F, {\Z_2}(i))]}$$ where, for $v\in
S\setminus S_\infty$,  $2^{n_v^{(i)}}$ is the exact power of $2$
dividing $e_{v}^{(i)}$ 
\qed 
\end{cor}

As mentioned before, if the number field is totally imaginary hypothesis   ($\mathcal{H}_i$)  always holds. 
However  ($\mathcal{H}_i$)  does not necessarily hold in some extension  $L/F$. 
For instance, take  $F=\Q$  and  $L=\Q(\sqrt{3})$. Then,  $L$ has only one dyadic prime and its class group is trivial. 
As explained in the general case in the proof of  Proposition  \ref{description 2,i regularity},  for such a number field we have 
$$H^1( o_L^\prime, \Z_2(i))/2\cong  H^1( o_L^\prime, \Z/2(i)) \cong  U_L^{\prime }/U_L^{{\prime }^2}.$$ 
Besides,  the group  $U_L^{\prime }$  of  $2$-units is generated by $-1, 2-\sqrt{3}$  and  $\sqrt{3}-1$.  Therefore the cokernel of the  signature map  $\mathrm{sgn}_L : H^1(L, {\Z_2}(i))/2 \longrightarrow  ({\Z}/2)^{r_1(L)}$  is  $\Z/2$ whose  $G$-cohomology groups  $\hat{H}^n(G,\Z/2)$  are non trivial. 
When the extension  $L/F$ is  cyclic, we have the following genus formula independent of hypothesis   ($\mathcal{H}_i$). 

\begin{thm}\label{odd cyclic genus}
Let  $L/F$ be a  cyclic extension of number fields with Galois
group  $G$  and  $i\geq 2$ odd.  Assume that $L/F$ is unramified at infinite primes.
Then
$$\frac{|H^2_\M(o_L, {\Z}(i))_{G}|}{|H^2_\M(o_F, {\Z}(i))|}= 
\frac{\prod_{v \in S \setminus S_\infty}e_{v}^{(i)}}{[H^1_\M(F,
{\Z}(i)): H^{1, \mathcal{N}}_\M(F, {\Z}(i))]}$$ where for each
finite prime  $v$, $e_{v}^{(i)} :=gcd (e_v, q_v^i-1)$,  $e_v$ being
the ramification index of  $v$   in the extension  $L/F$.
\end{thm}

\textbf{Proof}.  As in the proof of Theorem \ref{motivic genus}, we have a commutative diagram
$$
\begin{array}{ccc}
     {H}_1(G, H^2_\M(o_L^S,\Z(i)))& \stackrel{\alpha} \longrightarrow & {H}_1(G, \oplus_{w \mid v \in S\setminus S_\infty}H^1_\M(k_w,\Z(i-1))) \\
   \downarrow  \wr &  & \downarrow \wr\\
    \hat{H}^0(G, H^1_\M(L,\Z(i))) & \stackrel{\beta_S} \longrightarrow & \hat{H}^0(G, \oplus_{w \mid v \in S\setminus S_\infty}H^1_\M(k_w,\Z(i)))
 \end{array}
$$
by Remark  \ref{kcyclic shifting}.  
The rest of the proof goes along the same lines as Theorem  \ref{motivic genus}. \qed

For relative quadratic extensions,   explicit genus formulae are
given in \cite{Kolster03} using totally positive elements in the
Tate kernel (See also theorem  \ref{genus+}, \S 4 and \S 6). In
\cite{Assim-Movahhedi04, Assim-Movahhedi12}, 
we explicitly determined the above norm index in the particular cases where  the 
Galois group  $G$ is a cyclic  $p$-group for a prime  $p$. 
In section 4, we give the necessary and
sufficient condition for the above norm index to be exactly
$\prod_{v \in S\setminus S_\infty}e_{v}^{(i)}$. This allows us to study
the descent and co-descent for the motivic tame kernels in an
arbitrary Galois extension $L/F$.

\section{General genus formulae for even $K$-groups}
In  this section, we are going to apply the results and methods of
the previous sections to give genus formulae for even $K$-groups of rings of integers of number fields. We keep the notation of
the preceding sections:    $L/F$ is a finite Galois extension of
number fields with Galois group  $G$   and   $S$  is the set of
infinite primes and those which ramify in  $L/F$. For each finite
prime  $v$, $e_{v}^{(i)} :=gcd (e_v, q_v^i-1)$,  $e_v$ being the
ramification index of  $v$   in the extension  $L/F$. Finally $r$ is
the number of infinite primes of  $F$  which ramify  in the extension  $L/F$.

The relation  between $K$-theory and motivic cohomology is provided
by motivic Chern characters \cite{Pushin04}
 $$ch_{i,k}^\M : K_{2i-k} (F)   \longrightarrow  H^k_\M (F, {\Z}(i))$$
which, once tensored with  $\Z_p$, induce $p$-adic Chern characters
\cite[Chapter III]{Levine98}.  The odd part of the following theorem
is a consequence of  the motivic Bloch-Kato conjecture proved by Voevodsky
\cite{Voevodsky11}. The $2$-primary information  has been obtained by Kahn (partly) and
by Rognes-Weibel \cite{Kahn97, Weibel00} based on Voevodsky's proof
of the Milnor conjecture.

\begin{thm}\label{RW}
Let $S$ be a (finite) set of places of $F$ containing the places at
infinity and $o_F^S$ the ring of  $S$-integers of $F$.  Let $r_1$ be the
number of real places. Then for all $i\geq 2$ and $k=1,2$,   the
Chern characters
 $$ch_{i,k}^\M : K_{2i-k} (o_F^S)   \longrightarrow  H^k_\M (o_F^S, {\Z}(i))$$
are
\begin{enumerate}
\item[(i)] isomorphisms for $2i-k\equiv 0, 1, 2, 7
\pmod 8$
 \item[(ii)]  surjective with kernel isomorphic to
$(\Z/2)^{r_1}$ for   $2i-k\equiv 3 \pmod 8$
 \item [(iii)] injective with cokernel isomorphic to $(\Z/2)^{r_1}$ for $2i-k\equiv 6
\pmod 8.$
\end{enumerate}
If $i\equiv 3 \pmod 4$, there is an exact sequence
\begin{equation*}
0\rightarrow  K_{2i-1} (o_F^S)\rightarrow H^1_\M (o_F^S, {\Z}(i))
\rightarrow  ( {\Z}/2)^{r_1} 
\rightarrow K_{2i-2} (o_F^S) \rightarrow H^2_\M (o_F^S, {\Z}(i))
\rightarrow 0
\end{equation*}
where the middle map
$$H^1_\M (o_F^S, {\Z}(i)) \rightarrow  ( {\Z}/2)^{r_1}$$
factors through $H^1_\M (o_F^S, {\Z}(i))/2$ and induces the
signature map (\ref{sgn}).
\end{thm}

For each finite prime  $v$, we introduce a new map  $\varphi^{\prime}_v$ making
the following square commutative
$$
\begin{array}{cccccc}
 K_{2i-1} F & \stackrel{\varphi^{\prime}_v} \longrightarrow &   K_{2i-1} k_v\\
    \downarrow ch_{i,1}^\M &   &   \downarrow \wr      \\
H^1_\M(F, \Z(i)) & \stackrel{\varphi_v} \longrightarrow &
H^1_\M(k_v, \Z(i))
  \end{array}
$$
where  $\varphi_v$  has been introduced in section \ref{Genus
formula for motivic tame kernels}. The map   $\varphi^{\prime}_v$ will play a
similar role to  $\varphi_v$  in the context of  $K$-groups.

Consider the commutative diagram
$$
\begin{array}{cccccc}
 K_{2i-1} L  & \stackrel{\oplus_{w \mid v} \varphi^{\prime}_w} \longrightarrow & \oplus_{w \mid v} K_{2i-1} k_w \\
    \downarrow N_G &  &   \downarrow   N_G   \\
 K_{2i-1} F & \stackrel{\varphi^{\prime}_v} \longrightarrow &   K_{2i-1} k_v
  \end{array}
$$
Summing the cokernels of the vertical maps over all  $v \in S
\setminus S_\infty$, we have a canonical map
$$\beta_S^\prime : K_{2i-1} F /N_G (K_{2i-1}L )  \longrightarrow
\oplus_{v\in S \setminus S_\infty}  K_{2i-1} k_v /N_{G_v}
(K_{2i-1}k_w )$$ playing  a similar role to  $\beta_S$  defined in
section 2. Here, for each  $v$  we fix a $w$  in  $L$  above  $v$.
\textbf{\begin{defi} We define the normic subgroup  $K_{2i-1}^{\mathcal{N}}F$  to consist of
elements $a \in K_{2i-1} F$ whose image by $\varphi^{\prime}_v$  belongs to
 $N_{G_v}
(K_{2i-1}k_w )$  for all $v \in S \setminus S_\infty$ and $w\mid v$.
In other words, $K_{2i-1}^{ \mathcal{N}}F$ is the kernel of the map
$$K_{2i-1}F \longrightarrow
\oplus_{v\in S\setminus S_\infty}
 K_{2i-1} k_v /N_{G_v}
(K_{2i-1}k_w )$$ induced by the maps $\varphi^{\prime}_v$.
\end{defi}}

The transfer map
$$ \mathrm{Tr}_i  : K_{2i-2} (o_L^S)_G  {\rightarrow} K_{2i-2} (o_F^S)$$
realises isomorphisms on the odd  $p$-primary parts by Theorems
\ref{galois descent and codescent}(i) and \ref{RW}   provided that
$L/F$  is unramified outside  $S$. We are interested in the kernel
and cokernel of the transfer map at the level of the ring of
integers:
$$ \mathrm{Tr}_i  : K_{2i-2} (o_L)_G  {\rightarrow} K_{2i-2} (o_F)$$
In this section we prove genus formulae for even $K$-groups
combining the above Theorem \ref{RW} with Theorem \ref {motivic
genus} and Corollary \ref{odd 2-primary genus} from the previous
section.

\begin{thm}\label{genus K}
 Let $L/F$ be a finite Galois extension of number fields with Galois
group  $G$  and  $i\geq 2$ an integer. When  $i$  is odd, we assume hypothesis  ($\mathcal{H}_i$) holds. Then
$$\frac{|K_{2i-2} (o_L)_{G}|}{|K_{2i-2} (o_F)|}=  2^{\alpha_i}
\frac{\prod_{v \nmid \infty}e_{v}^{(i)}}{[K_{2i-1} F:
K_{2i-1}^{ \mathcal{N}}F]}$$ where
\begin{enumerate}
\item[(i)]  $s_i\leq \alpha_i \leq r$, if $2i-2\equiv 0 \pmod 8$;
\item[(ii)] $\alpha_i=-r$ if  $2i-2\equiv 2 \pmod 8$;
\item[(iii)]$\alpha_i=0$ if $2i-2\equiv 4  \; \mbox{or} \; 6 \pmod 8$.
\end{enumerate}

\end{thm}

%\item[(i)]$\alpha_i=\nu_i$ if  $2i-2\equiv
% 0 \pmod 8$
%\item[(iii)]$\alpha_i=1$ if $2i-2\equiv
% 4 \pmod 8$

%The composite map
%$$K_{2i-1}F \stackrel{ch_{i,1}^\M}{\longrightarrow} H^1_\M(F, {\Z}(i))
%\stackrel{\varphi_v}\rightarrow H^1_\M(k_v, {\Z}(i))$$

%for each $p$ a map
%$$H^1_\M(F, {\Z}(i)) \otimes \Z_p \rightarrow H^1_\M(k_w, {\Z}(i)) \otimes \Z_p.$$
%
%
%Let us first remark that for each finite place $v$, we have a canonical map
%$$K_{2i-1}F {\rightarrow} K_{2i-1}k_v.$$
%Indeed,  for each prime $p$, let  $p^m$ be the exact power of $p$ dividing the order $q_v^i-1$  of the $K$-group $K_{2i-1}k_v.$
%By (13), we have a map
%$$H^1_\M(F, {\Z}(i))/p^m{\rightarrow} H^1_\M (k_v, {\Z}(i))/p^m $$
%and hence a map
%$$H^1_\M(F, {\Z}(i))/(q_v^i-1){\rightarrow} H^1_\M (k_v, {\Z}(i))/(q_v^i-1)  \cong  H^1_\M(k_v, {\Z}(i)).$$
%Further the Chern character induces a map
%$$ K_{2i-1}F/(q_v^i-1) {\rightarrow} H^1 (F, {\Z_p}(i))/(q_v^i-1)$$
%
%The desired map is obtained as the composition
%$$ K_{2i-1}F {\rightarrow} K_{2i-1}F/(q_v^i-1) {\rightarrow} H^1 (F, {\Z_p}(i))/(q_v^i-1) {\rightarrow}
%K_{2i-1}k_v .$$

\proof
%Notice that since the \'etale Chern characters are isomorphisms for odd primes, we mainly focus our attention on   the $2$-primary part.

%\item Assume first that $2i-2\equiv 0
%\pmod 8$. The motivic Chern characters
%are isomorphisms and   a genus formula  for $K_{2i-2} (o_L)$ results immediately from
%Theorem \ref{motivic genus}. 
When $2i-2\equiv 0  \pmod 8$,  the motivic Chern characters 
$ch_{i,k}^\M : K_{2i-k} (o_F^S)   \longrightarrow  H^k_\M (o_F^S, {\Z}(i))$ 
from Theorem \ref{RW}  are isomorphisms for  $k=1, 2$ and the result
follows from Theorem \ref{motivic genus}.

When  $2i-2\equiv  2 \pmod 8$,  Theorem \ref{RW} shows that  $K_{2i-2}o_F \cong H^2_\M(o_F, {\Z}(i))$  and that we
have a commutative diagram
$$
\begin{array}{cccccccccccc}
   &  &   0 & \rightarrow
&(\Z/2)^{r_1(L)} & {\rightarrow} & K_{2i-1} L & \rightarrow &
H^1_\M(L, {\Z}(i))  & \fd & 0\\
           &  & & &   \downarrow &    &   \downarrow &             &   \downarrow \\
     & &  0 & \rightarrow & (\Z/2)^{r_1(F)} &
{\rightarrow} & K_{2i-1} F & \rightarrow & H^1_\M(F, {\Z}(i)) &
\rightarrow &  0
  \end{array}
  $$
where all the vertical maps are induced by the norm  $N_G$.  In
particular, we obtain an exact sequence
$$  (\Z/2)^r 
{\rightarrow}  K_{2i-1} F/N_G( K_{2i-1}L) \rightarrow  H^1_\M(F,
{\Z}(i))/N_G  (H^1_\M(L, {\Z}(i)) \rightarrow   0.$$ Hence a
commutatif diagram
$$
\begin{array}{cccccccccccc}
   &  &    &
& & & (K_{2i-1} F ) /N_G (K_{2i-1}L )&
 \twoheadrightarrow &
H^1_\M(F, {\Z}(i))/N_G  H^1_\M(L, {\Z}(i)) &  & \\
           &  & & &    &    &  \beta_S^\prime   \downarrow &             &  \beta_S \downarrow \\
     & &   &  & &
 & \oplus_{v\in S \setminus
S_\infty}  (K_{2i-1} k_v )/N_{G_v} (K_{2i-1}k_w )& \cong &
\oplus_{v\in S \setminus S_\infty} H^1_\M(k_w, {\Z}(i))/N_{G_v}H^1_\M(k_w,
{\Z}(i)) &  &
  \end{array}
  $$
showing that
$$ \mathrm{\Image}(\beta_S^\prime) \cong \mathrm{\Image}(\beta_S).$$
In particular, these two images have the same order:
$$
[K_{2i-1} F: K_{2i-1} ^{\mathcal{N}} F ] = [H^1_\M(F, {\Z}(i)):
H^{1,\mathcal{N}}_{\M}  (F, {\Z}(i))]. 
$$
 Again, the genus formula in this case results from Theorem \ref{motivic genus}.
% then $K_{2i-2} (o_L) \otimes \Z_2 \cong H^2(o_L[1/2], \Z_2(i))$

%
%If $2i-2\equiv 0, 2
%\pmod 8$,  the motivic Chern characters
%are isomorphisms and   a genus formula  for $K_{2i-2} (o_L)$ results immediately from
%Theorem \ref{motivic genus} :
%
%\begin{thm}
% Let $i\geq 2$ such that $2i-2\equiv 0, 2
%\pmod 8$. Then
%$$\frac{|K_{2i-2} (o_L)_{G}|}{|K_{2i-2} (o_F)|}= 2 ^{\nu_i}
%\frac{\prod_{v \in S \setminus S_\infty}e_{v}^{(i)}}{[H^1_\M(F,
%{\Z}(i))/n: H^1_\M(F, {\Z}(i))/n \cap (\bigcap_{v\in S \setminus
%S_\infty} N_{G_v}(H^1_\M(k_w, {\Z}/n(i)))]}$$ where   $\nu_i =-r$
%for $i$ even and $s_i\leq \nu_i \leq r$ for $i$ odd.
%\end{thm}

Assume now that $2i-2\equiv 4  \pmod 8$.  This is the case where
neither  $ch_{i,1}^\M$  nor $ch_{i,2}^\M$  is an isomorphism in
general.  By Theorem \ref{RW}, we have a commutative diagram
$$
\begin{array}{cccccccccccc}
&  &    & &0 &  & & &
 & &\\
           &  & & &  \downarrow  &    &   &             &   \\
   &  &    &
&((\Z/2)^{\delta_i(L)} )_G & {\rightarrow} & (K_{2i-2} o_L)_G&
 \rightarrow & H^2_\M(o_L, {\Z}(i))_G & \fd & 0\\

&  & & &  \downarrow  &    &  \mathrm{Tr}_i  \downarrow &             &  \mathrm{tr}_i \downarrow \\

&  & 0&\rightarrow &(\Z/2)^{\delta_i(F)}& \rightarrow & K_{2i-2} o_F & \rightarrow& H^2_\M(o_L, {\Z}(i))&\rightarrow &  0\\

 &  & & &  \downarrow  &    &   &             & \downarrow \\

&  & & & (\Z/2)^{s_i(L/F)}&  &  & & 0& &  \\

 &  & & &  \downarrow  &    &   &             &  \\

&  & & & 0 &  &  & & & &  \\
 \end{array}
  $$
where the right hand vertical map is surjective by Proposition \ref{beginning genus}. 
The left hand column is exact since  $\hat{H}^{-1}(G, (\Z/2)^{\delta_i(L)})=0$ by hypothesis  ($\mathcal{H}_i$).   
Applying Snake lemma, we get
\begin{equation} \label{comparing K4 genus motivic genus}
\frac{|K_{2i-2} (o_L)_{G}|}{|K_{2i-2} (o_F)|}= 2^{-s_i(L/F)}
\frac{|H^2_\M (o_L, {\Z}(i))_G |}{|H^2_\M (o_F, {\Z}(i))|} 
\end{equation}
for  $2i-2\equiv 4  \pmod 8$ under $(\mathcal{H}_i)$.

We  focus on the $2$-primary part since for the $p$-primary parts
with  $p$-odd the \'etale Chern characters realise isomorphisms
between $K_{2i-k}o_F \otimes {\mathbf{Z}}_{p}$  and   $H^{k}_\M
(o_F, \mathbf{Z}(i))\otimes {\mathbf{Z}}_{p}$  for  $k=1,2$. The
$2$-primary part of $\frac{|H^2_\M(o_L, {\Z}(i))_G |}{|H^2_\M(o_F,
{\Z}(i))|}$ has been computed in Corollary \ref{odd 2-primary
genus}. We shall prove that the two norm indices  
%\begin{equation} \label{comparing K5 norm + norm}
%[K_{2i-1} F: K_{2i-1} ^{\mathcal{N}} F ]  \; {_{2} }{\sim} \;
%[H^1_+(F, {\Z_2}(i)) : H^{1, {\mathcal{N}}}_+(F, {\Z_2}(i)) ]
%\end{equation} 
$[K_{2i-1} F: K_{2i-1} ^{\mathcal{N}} F ]$  and  
$[H^1_+(F, {\Z_2}(i)) : H^{1, {\mathcal{N}}}_+(F, {\Z_2}(i)) ]$
have the same  $2$-primary parts. 
Theorem \ref{RW} and the exact sequence (\ref{kolster2}) yield an
exact sequence
$$  0\rightarrow \oplus _{ v\, complex} \Z_2(i)\rightarrow H^1_+(L,
\Z_2(i)) \rightarrow  K_{2i-1} L \otimes \Z_2 \rightarrow 0.$$ We
then have a commutative diagram
$$
\begin{array}{cccccccccccc}
   &  &  0  &\rightarrow
&\oplus _{ w\mid v\, complex} \Z_2(i) & {\rightarrow} & H^1_+(L,
\Z_2(i)) &
 \rightarrow &
 K_{2i-1} L \otimes \Z_2 & \fd & 0\\
           &  & & &  \downarrow  &    &    \downarrow &             &   \downarrow \\
         &  &  0  &\rightarrow
&\oplus _{ v\, complex} \Z_2(i) & {\rightarrow} & H^1_+(F, \Z_2(i))
&
 \rightarrow &
 K_{2i-1} F \otimes \Z_2 & \fd & 0
  \end{array}
  $$
where the vertical maps are the norm maps  $N_G$. It follows that
$$ H^1_+(F, {\Z_2}(i)) /N_G H^1_+(L ,{\Z_2}(i)) \cong ( K_{2i-1} F \otimes \Z_2)/N_G( K_{2i-1} L \otimes \Z_2).$$
The  commutative diagram
$$
\begin{array}{cccccccccccc}
   &  &    &
&H^1_+(F, {\Z_2}(i)) /N_G H^1_+(L, {\Z_2}(i) & \cong & ( K_{2i-1} F
\otimes \Z_2)/N_G( K_{2i-1} L \otimes \Z_2) &
 &
  & & \\
           &  & & & \beta_S^+  \downarrow  &    &  \beta_S^\prime  \downarrow   &             &   \\
         &  &    & &
\oplus_{v\in S \setminus S_\infty} H^1(k_v,
{\Z_2}(i)/N_{G_v}H^1_+(k_w, {\Z_2}(i) & {\cong} & \oplus_{v\in S
\setminus S_\infty} (K_{2i-1} k_v \otimes \Z_2)/N_{G_v}(K_{2i-1} k_w
\otimes \Z_2) &
 &
 &  &
  \end{array}
  $$
gives the desired equality between norm indices.

Assume finally that $2i-2\equiv 6 \pmod 8$.
By  Theorem \ref{RW}(iii), we have a commutative diagram
$$
\begin{array}{cccccccccccc}
          &  & & &   &    &   &             &    0 \\
          &  & & &   &    &   &             &    \downarrow \\
   &  &  H_1(G,\oplus_{w \, real} \Z/2) =0  &{\rightarrow}
&K_{2i-2} (o_L)_G & {\rightarrow} & H^2_\M(o_L, {\Z}(i))_G &
\rightarrow &
(\oplus_{w \, real} \Z/2)_G & \fd & 0\\
           &  & & &  \mathrm{Tr}_i \downarrow &    &   \mathrm{tr}_i \downarrow &             &    \downarrow \\
     & &  0 & \rightarrow & K_{2i-2} (o_F) &
{\rightarrow} & H^2_\M(o_F, {\Z}(i))& \rightarrow & \oplus_{v  \, real} \Z/2 & \rightarrow &  0 \\
           &  & & &   &    &   \downarrow &             &    \downarrow \\
     & &   &  &  &
 & (\Z/2)^r & = & (\Z/2 )^r &  \\
          &  & & &   &    &   \downarrow &             &    \downarrow \\
     & &   &  &  &
 &0 & & 0&
  \end{array}
$$
The exactness of the middle vertical sequence comes from Proposition
\ref{beginning genus}. Hence the  transfer map  $ \mathrm{Tr}_i $ is
surjective and  has the same kernel as the map  $ \mathrm{tr}_i$.
Accordingly
\begin{equation}
\label{K6 motivic} \frac{|K_{2i-2} (o_L)_{G}|}{|K_{2i-2} (o_F)|}=
2^{r} \frac{|H^2_\M(o_L, {\Z}(i))_G| }{|H^2_\M(o_F, {\Z}(i))| }
\end{equation}
for $2i-2\equiv 6  \pmod 8$.

To conclude, we use again the genus formula in Theorem \ref{motivic
genus}  and the fact that  $K_{2i-1} (F)  \cong H^1_\M(F, \Z(i))$.
\qed

\section{Galois co-descent} \label{S galois codescent}
%Let  $L$ be a finite Galois extension  of $F$  with Galois
%group  $G$    unramified outside a set of primes $S$ of  $F$.
In this section, we will carry out an arithmetic interpretation concerning Galois co-descent
for motivic tame kernels  $H^2_\M(o_F, {\Z}(i))$. 

We first give a necessary and sufficient condition under which the map
$$\beta_S : \hat{ H }^0 (G, H^1_\M (L,\Z(i))) \rightarrow   \oplus_{v\in S\setminus S_\infty} \hat{ H }^0 (G_v, H^1_\M(k_w,\Z(i))) $$
introduced in section \ref{Genus formula for motivic tame kernels}  is surjective. 
For each prime $p$, let
$$T_p= \{v \in S \setminus S_\infty / \;  p\mid e_v\}$$
and
$$ \mathrm {loc}_{T_p} : H^1_\et( F,\Z_p(i)) \rightarrow  \oplus_{v\in T_p} H^1_\et(k_v,\Z_p(i))$$ 
be the direct sum   $\oplus_{v \in T_p} \varphi_v$ tensored with  $\Z_p$. 
For the definition of the map  $\varphi_v$  see the \S 2. 
%Introduce the integer
%$$e:=lcm \{e_v / \; v\in S\setminus S_\infty\}.$$

%{\bf Introduire
%$$\beta_{S,p} : \hat{ H }^0 (G,  H^1_\et( L,\Z_p(i))) \rightarrow  \oplus_{v\in S\setminus S_\infty} \hat{ H }^0 (G_v,H^1_\et(k_w,\Z_p(i))).$$
%for  each prime  $p$ et expliquer   $\beta_S = \oplus_{p\mid e}\beta_{S,p}$.}
%
%
%With this notation:

\begin{lem}\label{betas loc}
The map $\beta_S$  is surjective precisely when all the localization maps $ \mathrm {loc}_{T_p}$  are
surjective for all primes $p$ (dividing the ramification index $e_v$  for some   $v\in S \setminus S_\infty$). 
%The following two conditions are equivalent \\
%(i)  the map $\beta_S$  is surjective; \\
%(ii) the localization map $ \mathrm {loc}_{T_p}$  is
%surjective for each prime $p$ dividing  $e$.
\end{lem}

{\bf Proof.} 
For each prime  $p$ and each finite prime $v$, we have a commutative diagram 
%\begin{equation*}
% \begin{array}{cccccccc}
%H^1_\et(L,\Z_p(i)) & \stackrel{N_G} \longrightarrow &  H^1_\et(F,\Z_p(i)) &  \longrightarrow & \hat{ H }^0 (G,  H^1_\et( L,\Z_p(i))) &\longrightarrow & 0 \\
%\downarrow &  & \downarrow {\mathrm {loc}_{T_p} &  & \downarrow &&\\
%\oplus_{v \in {T_p}\oplus_{w \mid v} H^1_\et(k_w,\Z_p(i)) & \stackrel{N_G}  \longrightarrow  &   \oplus_{v \in {T_p} H^1_\et(k_v,\Z_p(i)) &  \longrightarrow & \oplus_{v\in T_p} \hat{ H }^0 (G_v,H^1_\et(k_w,\Z_p(i))) &  \longrightarrow & 0
%\end{array}
% \end{equation*}
 %
%Let  $p$ be a prime dividing the ramification index  $e_v$  for some finite prime  $v$.
%
\begin{equation*}
 \begin{array}{cccccc}
H^1_\et(L,\Z_p(i)) &  \longrightarrow &   \oplus_{w \mid v} H^1_\et(k_w,\Z_p(i)) & &\\
\downarrow N_G&  & \downarrow \oplus N_{L_w/k_v}  & & \\
H^1_\et(F,\Z_p(i)) &  \stackrel{\varphi_v \otimes \Z_p} \longrightarrow &   H^1_\et(k_v,\Z_p(i))  & & \\
        \downarrow &  & \downarrow  & & \\
\hat{ H }^0 (G,  H^1_\et( L,\Z_p(i))) &  \longrightarrow &  \hat{ H }^0 (G, \oplus_{w \mid v} H^1_\et(k_w,\Z_p(i)))  & \cong &  \hat{ H }^0 (G_v,H^1_\et(k_w,\Z_p(i))) \\
\end{array}
 \end{equation*}
%\begin{equation*}
% \begin{array}{ccc}
%     H^1_\et(L,\Z_p(i)) &  \longrightarrow &   \oplus_{v \in {T_p}\oplus_{w \mid v} H^1_\et(k_w,\Z_p(i))\\
%   \downarrow N_G&  & \downarrow \oplus N_G\\
%     H^1_\et(F,\Z_p(i)) &  \longrightarrow &   \oplus_{v \in {T_p} H^1_\et(k_v,\Z_p(i))\\
%        \downarrow &  & \downarrow \\
%%H^1_\et(F,\Z_p(i))/N_G(H^1_\et(L,\Z_p(i))) &  \longrightarrow & \oplus_{v\in {T_p} H^1_\et(k_v,\Z_p(i))/N_G(\oplus_{w \mid v}H^1_\et(k_w,\Z_p(i)))\\
%%   \downarrow \wr&  & \downarrow \wr\\
%\hat{ H }^0 (G,  H^1_\et( L,\Z_p(i))) &  \longrightarrow & \oplus_{v\in T_p} \hat{ H }^0 (G, \oplus_{w \mid v} H^1_\et(k_w,\Z_p(i))) \\
%&  & \downarrow \wr\\
%&  \longrightarrow & \oplus_{v\in T_p} \hat{ H }^0 (G_v,H^1_\et(k_w,\Z_p(i))) \\
%\end{array}
% \end{equation*}
%
%
%
%Let  $p$ be a prime dividing  $e$.
%\begin{equation*}
% \begin{array}{ccc}
%     H^1_\et(L,\Z_p(i)) &  \longrightarrow &   \oplus_{v \in {T_p}\oplus_{w \mid v} H^1_\et(k_w,\Z_p(i))\\
%   \downarrow N_G&  & \downarrow \oplus N_G\\
%     H^1_\et(F,\Z_p(i)) &\stackrel{\mathrm {loc}_{T_p}  \longrightarrow &   \oplus_{v \in {T_p} H^1_\et(k_v,\Z_p(i))\\
%\end{array}
% \end{equation*}
where the vertical maps are exact since  $H^1_\et( L,\Z_p(i))$ and $ H^1_\et(k_w,\Z_p(i))$ 
satisfy Galois descent  (Theorem \ref{galois descent and codescent} (i) and  the isomorphism (\ref{residual descent})). 

Suppose first that all the maps  $\mathrm {loc}_{T_p}$  are surjective. 
%{\bf If  $p$ does not divide $e_v$, as seen in the proof of Lemma
%\ref{lemma evi},  the Galois cohomology group $\hat{ H }^0 (G_v,
%H^1_\et(k_w,\Z_p(i)))$ vanishes. So we have to concentrate on those primes $p$ which divide  $e_v$  for some finite prime  $v$ \\
%a-t-on besoin de le dire et le dire ici ?}. 
Then, by summing the above diagram over all  $v\in T_p$, we see that 
$$\beta_{T_p} : \hat{ H }^0 (G,  H^1_\et( L,\Z_p(i))) \rightarrow  \oplus_{v\in T_p} \hat{ H }^0 (G_v,H^1_\et(k_w,\Z_p(i)))$$ 
is also surjective. 
Hence the surjectivity of  $\beta_S$  which is, by definition, the direct sum  $\oplus_p\beta_{S,p}$. 

Conversely, suppose that  $ \beta_S$ is surjective. 
For each prime $v$
in $S\setminus S_\infty$, fix a prime $w$ of  $L$ above $v$. As in
the proof of Lemma \ref{lemma evi}, we have
$$ \hat{H}^0(G_v, H^1_\M(k_w,\Z(i))) \cong  H^1_\M(k_v,\Z(i))/e_v. $$ 
%
%Let  $I_v$  be the inertia group of $v$ in the extension $L_w/F_v$.
% By \cite[section 12, page 185 remark]{Quillen72}, the norm map
%  $$N_{G_v/I_v}: H^1_\M(k_w,\Z(i))  \rightarrow   H^1_\M(k_v,\Z(i))$$
%  is surjective and, hence, induces a surjective map
%  $$H^1_\M(k_w,\Z(i))/e_v  \rightarrow   H^1_\M(k_v,\Z(i))/e_v.$$
%Now,   $I_v$  acting trivially on  $H^1_\M(k_w,\Z(i))$, we have
%$H^1_\M(k_w,\Z(i))/e_v \cong  \hat{H}^0(I_v, H^1_\M(k_w,\Z(i)))$ so that
% $$
%  \begin{array}{lll}
%      \hat{H}^0(G_v, H^1_\M(k_w,\Z(i))) & \cong &  H^1_\M(k_v,\Z(i))/ N_{G_v/I_v}(N_{I_v} H^1_\M(k_w,\Z(i))) \\
%        & \cong & H^1_\M(k_v,\Z(i))/e_v N_{G_v/I_v}( H^1_\M(k_w,\Z(i)))\\
%        & \cong &  H^1_\M(k_v,\Z(i))/e_v.
%     \end{array}
%  $$
Therefore the following map   
   $$H^1_\M (F,\Z(i)) \rightarrow \oplus_{v\in S\setminus S_\infty} H^1_\M(k_v,\Z(i))/e_v ,$$
obtained naturally from the $\varphi_v,  v\in S\setminus S_\infty$,  is  surjective. Dividing the above map by $p$, we get the
surjectivity of the map 
   $$H^1_\M (F,\Z(i))/p \rightarrow \oplus_{v\in T_p} H^1_\M(k_v,\Z(i))/p $$
which in turn shows that of  $\mathrm {loc}_{T_p}$. \qed

Suppose first that  $i\geq 2$ is even. In this case
$\mathrm{coker}(\alpha) \cong \mathrm{coker}(\beta_S)$ (see the
discussion before Theorem \ref{motivic genus}). It follows by
Proposition \ref{beginning genus} that the sequence
$$0 \rightarrow H^2_\M(o_L, {\Z}(i))_{G} {\rightarrow}
H^2_\M(o_F, {\Z}(i)) {\rightarrow} ({\Z}/2)^r \rightarrow 0$$ is
exact precisely when the equivalent conditions of the above lemma
hold.

If instead  $i \geq 3$ is odd, we assume that hypothesis $(\mathcal{H}_i)$  holds (or $L/F$ cyclic)  and that 
%and  that the partial signature
%map
%$$H^1_\M(F, {\Z}(i))/2  \rightarrow  ({\Z}/2)^r$$
%is trivial (i.e $r=s_i$). This happens for instance when 
our
extension  $L/F$ is unramified at infinity. Hence, once again
$\coker(\alpha) \cong \coker(\beta_S)$ (see the
discussion before Theorem \ref{motivic genus}) and by Proposition
\ref{beginning genus} 
%the sequence
%$$0 \rightarrow ({\Z}/2)^{s_i} \rightarrow 
$$H^2_\M(o_L, {\Z}(i))_{G} \, {\cong} \,
H^2_\M(o_F, {\Z}(i))  $$
%\rightarrow 0$$ 
%is exact 
precisely when the
equivalent conditions of the above lemma hold.

For a prime $p$, let  $E=F(\mu_p)$  and  $\Delta =
\mathrm{Gal}(E/F)$. Let  $v$  be a finite non-$p$-adic prime of $F$
which splits in  $E$. Let  $v^\prime$  be any prime of  $E$ above
$v$. If $M/E$  is a Galois extension unramified at  $v^\prime$, then
the  Frobenius automorphism in $M/E$ does not depend on the choice
of $v^\prime$ above  $v$  and will simply be denoted by
$\sigma_v(M/E)$.

We now recall the definition of the \'etale Tate kernel  $D_F^{(i)}$  introduced in  \cite{Kolster91}. 
Let  $E:=F(\mu_p)$  and  $\Delta :=\Gal(E/F)$. 
For $i \geq 2$,  the exact sequence
$$0 \rightarrow {\bf Z}_p(i) \rightarrow {\bf Z}_p(i) \rightarrow
{\bf Z}/p{\bf Z}(i) \rightarrow 0$$ induces, by cohomology, the
following commutative diagram
$$\begin{array}{ccccccccccc}
0 & \rightarrow &H^1(F,{\Z}_p(i)) /p & \rightarrow &
H^1(F,{\Z}/p(i)) & \rightarrow &
_{p}H^2(F,{\Z}_p(i)) &\rightarrow & 0 \\
& & {\downarrow}\wr & & {\downarrow}\wr & & {\downarrow}\wr & \\
0 & \rightarrow & (H^1(E,{\Z}_p(i)) /p )^\Delta & \rightarrow &
E^{\bullet}/E^{\bullet p}(i-1)^{\Delta} & \rightarrow &
(_{p}H^2(E,{\Z}_p(i)))^\Delta &\rightarrow & 0.
\end{array}$$

This shows the existence of a subgroup  $D_F^{(i)}$  of
$E^{\bullet}$  containing  $E^{\bullet p}$ - the analogue of the
Tate kernel in the case where  $i=2$, $F\supset \mu_p$ -, such that
$$H^1_\et( F,\Z_p(i))/p   \simeq D_F^{(i)}/E^ {\bullet p}(i-1).$$
Let  $G_{S_p}(F)$  be the Galois group over $F$ of the maximal
extension of  $F$  unramified outside  $p$ and infinity. 
Since
%\cite{Soule79}   {\bf Milne pour le deuxime isomorphisme, arithmetic duality thm,  (donner la ref exacte : peut-tre KNF}
$$H^1_\et( F,\Z_p(i)) \cong H^1_\et(o_F^\prime,\Z_p(i)) \cong H^1(G_{S_p}(F),\Z_p(i))$$
we have
$$ D_F^{(i)}/E^ {\bullet p} (i-1) \hookrightarrow H^1(G_{S_p}(F),\Z/p(i)) \cong H^1(G_{S_p}(E),\Z/p(i))^\Delta .$$
%Following the methods of Greenberg \cite{Greenberg78}, one shows \cite{Assim-Movahhedi04,  Assim-Movahhedi12} that
%$$D_{F}^{(i)} = \{a\in E^{\bullet}  / a \otimes p^{-1} \in
%{\mathrm{Div}} ({\mathcal M}(i-1)^{G_\infty})\}$$
%where  ${\mathcal M}$  is the Kummer radical of the maximal abelian pro-$p$-extension of  $E_\infty$ unramified outside $p$ and  $G_\infty= Gal(E_\infty/F)$.
Therefore the elementary extension $E(\sqrt[p]{D_F^{(i)}})/E$ is
unramified outside $p$  and infinity.

This being said,  since  $H^1_\et(k_v,\Z_p(i))=0$ for any $p$-adic
prime $v$,  the localization map $\mathrm {loc}_{T_p}$ can be
written
$$\mathrm {loc}_{T_p} : H^1_\et( F,\Z_p(i)) \rightarrow  \oplus_{v\in T_p \setminus S_p} H^1_\et(k_v,\Z_p(i)).$$
Its surjectivity is equivalent to that of
$$ H^1_\et( F,\Z_p(i))/p \rightarrow  \oplus_{v\in T_p\setminus S_p} H^1_\et(k_v,\Z_p(i))/p .$$
For each  $v\in T_p\setminus S_p$, we have  $\mu_p \subset
F_v$ and
%( $p$  divise  $G_0/G_1 \hookrightarrow k_v^\bullet$  CL Serre)
$$H^1_\et(k_v,\Z_p(i))/p \cong  H^1(k_v,\Z/p(i)) \cong H^1(F_{v,\infty}/F_v,\Z/p(i))$$
where  $F_{v,\infty}:= F_v(\mu_{p^\infty})$  is the cyclotomic  $\Z_p$-extension of the local field  $F_v$ which coincides with the maximal
unramified pro-$p$-extension of  $F_v$. 
Therefore, the last map is the dual of the map
$$\varphi^{\prime}_{T_p} :  \prod_{v\in T_p \setminus S_p}
\mathrm{Gal}(F_{v,\infty}/F_v)/p \longrightarrow
\mathrm{Gal}(E(\sqrt[p]{D_F^{(i)}})/E) $$ sending each term
$\mathrm{Gal}(F_{v,\infty}/F_v)/p$  onto the decomposition group of
the prime  $v$  in the extension  $E(\sqrt[p]{D_F^{(i)}})/E$.

Hence, the surjectivity of  $\mathrm {loc}_{T_p}$  is equivalent  to the
injectivity of  $\varphi^{\prime}_{T_p}$. Summarizing we obtain the
following Galois co-descent criterion for motivic tame kernels: 

\begin{thm}\label{motivic codescent thm}
Let   $L/F$ be a finite Galois extension  of number fields  with
Galois  group  $G$ and $i\geq 2$. 
  Assume that $L/F$ is unramified at infinity. If $i$ is odd, we assume either hypothesis $(\mathcal{H}_i)$ holds or $L/F$  is cyclic. Then,   the natural map tr$_i : H^2_\M(o_L,
{\Z}(i))_{G} \rightarrow H^2_\M(o_F, {\Z}(i))$ is surjective and the
following two conditions are equivalent
\begin{enumerate}
\item[(i)] the map tr$_i : H^2_\M(o_L, {\Z}(i))_{G}
\rightarrow H^2_\M(o_F, {\Z}(i))$ is an isomorphism
 \item [(ii)] for every prime $p$  dividing  the ramification index $e_v$  for some finite prime  $v$, the Frobenius automorphisms  $\sigma_v(E(\sqrt[p]{D_F^{(i)}})/E)$, $v\in T_p\setminus S_p$,  are linearly independent in the $\mathbf{F}_p$-vector space  $\mathrm{Gal}(E(\sqrt[p]{D_F^{(i)}})/E)$. 
 \end{enumerate}
Here  $T_p$  stands for the set of finite primes  $v$  of  $F$
such that  $p\mid e_v$. \qed
\end{thm}

\begin{rem}
The hypothesis "$L/F$ is unramified at infinity" in the above
theorem means that  $r=0$. When $i$  is even, it is a necessary
condition for the map tr$_i$  to be an isomorphism (see Proposition
\ref{beginning genus}). Hence, for  $i$  even,   tr$_i$ is an
isomorphism precisely when $r=0$  and Condition (ii) of the above
Theorem holds.
\end{rem}

 Now, let
${\widetilde{E}_1}$ be the compositum of the first layers of all
$\Z_p$-extensions of $E$. Denote by  $A_E$  its Kummer radical:
${\widetilde{E}_1} = E(\sqrt[p]{A_E})$. The Galois group  $\Delta$
acts on $\mathrm{Gal}({\widetilde{E}_1}/E)$ by conjugation. For
$i\in \Z$, let ${\widetilde{E}_1}^{(i)}$ be the subfield of
${\widetilde{E}_1}$ corresponding to
$\mathrm{Gal}({\widetilde{E}_1}/E)(-i)_\Delta$ and
$A_E^{[1-i]}:=A_E{(i-1)^\Delta}$  its Kummer dual:
$${\widetilde{E}_1}^{(i)} = E(\sqrt[p]{A_E^{[1-i]}}).$$

%$$ \mathrm{Gal}({\widetilde{E}_1}/E)(-i)_\Delta \simeq \mathrm{Gal}({\widetilde{E}_1}^{(i)}/E)\simeq
%\mathrm{Gal}(E(\sqrt[p]{A_E^{[1-i]}})/E).$$
%%$$\mathrm{Gal}(E(\sqrt[p]{A_F^{(i)}})/E)\simeq \mathrm{Gal}({\widetilde{E}_1}/E)(-i)_\Delta.$$
%$$\mathrm{Gal}({\widetilde{E}_1}^{(i)}/E)\simeq \mathrm{Gal}({\widetilde{E}_1}/E)(-i)_\Delta.$$
Let $T$ be a set of primes of $F$ containing the set $S_p$  of
primes of $F$ above $p$ and infinite primes. We assume that the
absolute norm  $Nv \equiv 1 \pmod p$ for all non-$p$-adic primes of
$T$ since  such a condition is necessary for  $e_v$  to be divisible
by  $p$.

Let  $i \in \Z$. The  set $T$  is called $i$-primitive for  $(F,p)$
if the Frobenius automorphisms
$\sigma_v({\widetilde{E}_1}^{(i)}/E)$, $v\in T\setminus S_p$,
generate an $\mathbf{F}_p$-subspace of
$\mathrm{Gal}({\widetilde{E}_1}^{(i)}/E)$ of dimension the
cardinality of $T\setminus S_p$. The notion of $i$-primitive sets
introduced in \cite{Assim95, Assim-Movahhedi12} is the twisted
version of the notion of primitive sets  \cite{Gras-Jaulent89,
Movahhedi88, Movahhedi90,  Movahhedi-Nguyen90}. When $i\equiv 0
\pmod \Delta$, $i$-primitive sets are exactly primitive sets. 
Here, we will use the following more suggestive definition of primitivity (see also \cite {Hutchinson05}):

%More generally, one can replace  $A_E^{[1-i]}$  by any subgroup  $D$  of  $E^\bullet$
%%containing  $E^{\bullet p}$
%such that  $E(\sqrt[p]{D})/E$ is unramified outside  $p$-adic primes as suggested in \cite {Hutchinson05}:

\begin{defi}\label{defi primitive} 
Let  $T$  be a set of  primes of $F$  containing $S_p$  and   $D$  a subgroup of  $E^\bullet$  such that  $E(\sqrt[p]{D})/E$  is unramified at all  $v\in T\setminus S_p$. Then  $T$  is called $D$-primitive for $(F,p)$ if the Frobenius automorphisms
$\sigma_v(E(\sqrt[p]{D})/E)$, $v\in T\setminus S_p$, 
are linearly independent in the  $\mathbf{F}_p$-vector space  $\mathrm{Gal}(E(\sqrt[p]{D})/E)$. 
\end{defi}

With this definition, the second condition of Theorem \ref{motivic
codescent thm}  means that for each prime $p$ the set  $T_p$
is $D_F^{(i)}$-primitive for $(F,p)$. By \u{C}ebotarev density
theorem, there exists an infinite number of $D$-primitive sets as
soon as $E(\sqrt[p]{D}) \neq E$.

%{\bf ceci n'est pas juste car parfois  $S_p$  est le seul ensemble primitif}

%In fact for such a number
%field, the Galois group of the maximal pro-$p$-extension of $F$,
%unramified ouside $p$, is a a free pro-$p$-group of rank $1 +r_2$.

%{\bf (citer kolster) }

Number fields such that the maximal pro-$p$-quotient of  $G_{S_p}(F)$  is a free pro-$p$-group are called $p$-rational. 
They were introduced in  \cite{Movahhedi88, Movahhedi90, Movahhedi-Nguyen90} to construct infinitely many examples of non-abelian extensions of  $\Q$  satisfying Leopoldt's conjecture at the prime $p$. 
More recently they have been used in order to construct continuous representations  
$\rho : G({\bar \Q}/\Q) \to  GL_n(\Z_p)$  
of the absolute Galois group of  $\Q$  with an open image for some $n \geq 2$  and  $p$ \cite{Greenberg16}. 
When  $p$  is odd, the number field  $F$  is  $p$-rational precisely when  
$ H^2_\et( o_F^\prime, {\Z/p})=0$.  
Number fields $F$ such that $H^2_\et( o_F^\prime,
{\Z/p}(i))=0$ are called $(p,i)$-regular  \cite{Assim95}. For  $p$
odd, the $(p,i)$-regularity is equivalent to the nullity of
$H^2_\et( o_F^\prime, {\Z_p}(i)).$ These number fields can be
considered  as a generalization of the cyclotomic fields
$\Q(\mu_p)$, $p$ regular. In the same way,
$(p,i)$-regular number fields are introduced to construct number
fields satisfying the "twisted" Leopoldt conjecture
%\cite{Greenberg78, Schneider79}:
$$(C_i): \; \; \,\;  H^2_\et( o_F^\prime, {\Q_p/\Z_p}(i))=0, \; \, i\not= 1.$$ 
As remarked in \cite{Kurihara92}, the nullity of the \'etale cohomology groups  $H^2_\et(\Z[\frac{1}{p}], {\Z_p}(i))$  for  all odd $i \geq 3$ is equivalent to Vandiver's conjecture for the prime $p$ (asserting the non divisibility by  $p$  of the class group of the maximal real subfield of  $\Q(\mu_p)$). For any given  odd $i$, Soul\'e proved that  $H^2_\et(\Z[\frac{1}{p}], {\Z_p}(i))=0$ for  $p$  greater than an effective bound depending on $i$  \cite{Soule99}.  The vanishing of  $K_4\Z=0$  is already known  (\cite{Rognes00}, see also \cite{Soule00}).

In \cite{Assim95} a going-up property for   $(p,i)$-regularity is
given for a  $p$-extension $L/F$  under the Leopoldt conjecture in
the cyclotomic tower $L(\mu_{p^\infty})$.  As an immediate consequence of  Theorem \ref{motivic
codescent thm}, we have an alternative proof of the same
property without any assumption on Leopoldt's conjecture. Recall that for  $(p,i)$-regular
number fields,  $i$-primitivity coincides with
$D_F^{(i)}$-primitivity  \cite[Prop. 1.2 and 1.4]{Assim-Movahhedi12}.

\begin{cor}\label{going up p,i regularity}
Let  $p$  be an odd prime number and let  $L/F$ be  a finite Galois
$p$-extension  of number fields  with Galois  group
  $G$. Then  for $i \in \Z$, the following conditions are equivalent
\begin{enumerate}
\item[(i)]  $H^2_\et(o_L^\prime, \Z_p(i))=0 ;$
\item[(ii)] $ H^2_\et( o_F^\prime, \Z_p(i))=0$ and the set of primes of $F$ which are in $S_p$ or ramify in   $L$ is $D_F^{(i)}$-primitive for $(F,p)$.
\end{enumerate}
\end{cor}

\textbf{Proof} It is clear that when $i\equiv j \mod  \mid \Delta
\mid$, the properties (i) and (ii) are exactly the same for $i$ and
$j$. So without loss of generality we can assume that $i\geq 2$  and
apply Theorem \ref{motivic codescent thm}. \qed

We also need a co-descent  criterion for the positive  cohomology groups $H^2_+ (o_L[1/2], {\Z_2}(i))$.
%We focus on the case where $L/F$ is a $2$-extension.
%If $i$ is even,  Proposition 2.8  and Theorem 2.4 show that
%$$H^2_+ (o_L[1/2], {\Z_2}(i))_G \cong H^2_+ (o_F[1/2], {\Z_2}(i))$$
% precisely when the Frobenius automorphisms
%$\sigma_v(F(\sqrt{D_F^{(i)}}/F)$, $v\in T_2\setminus S_2$,
%generate an $\mathbf{F}_2$-subspace of
%$\mathrm{Gal}(F(\sqrt{D_F^{(i)}})/F)$ of
% dimension the cardinality of $T_2\setminus S_2$.
%
%
%If $i$ is odd, 
Following Kolster \cite [Definition 2.4]{Kolster03},
we introduce the subgroup   $D_F^{+(i)}$  (the positive \'etale Tate
kernel)  to be the kernel of the signature map $\sgn$ (see Subsection 1.2)
restricted to  $D_F^{(i)}$ so that we have an exact sequence
\begin{equation}
 \label{sequence D+ and D}
0 \to D_F^{+(i)}/F^ {\bullet 2}  \to   D_F^{(i)}/F^ {\bullet 2}  \to (\Z/2)^{r_1(F)} \to (\Z/2)^{\delta_i(F)} \to 0.
\end{equation} 

Recall that, when $i$  is even,  Proposition \ref{sgn triviality} amounts to  $D_F^{+(i)} = D_F^{(i)}$. 

The same arguments as in the proof of Theorem 4.2 above show that
$$H^2_+ (o_L[1/2], {\Z_2}(i))_G \cong H^2_+ (o_F[1/2], {\Z_2}(i))$$
precisely when  the localization composite map
$$ H^1_+( F,\Z_2(i)) \rightarrow  H^1( F,\Z_2(i)) \stackrel{\mathrm {loc}_{T_2}} \longrightarrow  \oplus_{v\in T_2} H^1_\et(k_v,\Z_2(i))$$
is surjective, or equivalently the induced composite map
$$H^1_+( F,\Z_2(i)) /2  \rightarrow  H^1( F,\Z_2(i)) /2 \simeq D_F^{(i)}/F^ {\bullet 2}(i-1)  \to \oplus_{v\in T_2} H^1_\et(k_v,\Z_2(i))/2$$
is surjective. The image of the above left map being the subgroup
$D_F^{+(i)}/F^ {\bullet 2}$  of   $D_F^{(i)}/F^ {\bullet 2}$, we
have the following theorem which is proved in the same way as
Theorem \ref{motivic codescent thm}.

\begin{thm}\label{motivic+ codescent thm}
Let   $L/F$ be a  Galois $2$-extension  of number fields  with
Galois  group  $G$. Then, the surjective map $\mathrm{tr}_i ^+ :
H^2_+(o_L^\prime, {\Z_2}(i))_{G} \rightarrow H^2_+(o_F^\prime,
{\Z_2}(i))$  is an isomorphism precisely when $S_2 \cup T_2$  is
$D_F^{+(i)}$-primitive for  $(F,2)$.
%Here  $T_2$  stands for the set of finite primes  $v$  of  $F$  such that  $2\mid e_v$. \qed
\qed 
\end{thm}

\begin{cor}\label{application motivic+ codescent thm}
Let   $L/F$ be an imaginary $2$-extension  of number fields  with
    Galois  group  $G$ and $i\geq 2$ odd. Then the surjective map $\mathrm{tr}_i : H^2(o_L^\prime,
    {\Z}_2(i))_{G} \rightarrow H^2(o_F^\prime, {\Z}_2(i))$ is an
    isomorphism
    precisely when the signature map $\mathrm{sgn}_F : H^1(F, {\Z_2}(i))/2 \longrightarrow  ({\Z}/2)^{r_1(F)}$ is surjective and
    $S_2 \cup T_2$  is  $D_F^{+(i)}$-primitive for  $(F,2)$.
%Here  $T_2$  stands for the set of finite primes  $v$  of  $F$  such that  $2\mid e_v$. 
\end{cor}
%The above theorem leads immediately to the following corollaries :
%begin{cor}\label{going up 2,i regularity}
%Let  $L/F$ be  a finite Galois $2$-extension  of number fields  with
%Galois  group  $G$. Then  $H^2_+(o_L^\prime, \Z_2(i))=0$  precisely
%when $ H^2_+( o_F^\prime, \Z_2(i))=0$ and
%\begin{enumerate}
%\item[(i)]  if $i\geq 2$ is even, the set of primes of $F$ which are in $S_2$ or ramify in   $L$ is $D_F^{(i)}$-primitive for $(F,2)$;
%\item[(ii)] if $i\geq 2$ is odd, the set of primes of $F$ which are in $S_2$ or ramify in   $L$ is $D_F^{+(i)}$-primitive for $(F,2)$.
%\end{enumerate}
%\end{cor}

%To be complete
\textbf{Proof.} Since $L/F$ is totally imaginary,  hypothesis
($\mathcal{H}_i$) is satisfied and $s_i(L/F) = \delta_i(F)$. The exact
sequence (\ref{tr tr+}) shows that tr$_i$ is an isomrphism  exactly
when $s_i(L/F)=\delta_i(F) =0$ and tr$_i^+$ is an isomorphism. \qed

Let us mention the following arithmetic criterion for the nullity of the positive \'etale cohomology groups $H^2_+( o_F^\prime, \Z_2(i))$  whose proof follows the same ideas as those in the proof of  \cite[Prop.2.6]{Kolster03} where the case of real number fields is dealt with. See also \cite[Prop.2.2]{Ostvaer00}. 

\begin{proposition}\label{description 2,i regularity}
Let  $F$ be  a number field and $i\geq 2$. Then  $H^2_+(o_F^\prime, \Z_2(i))$  vanishes precisely when $F$ has only one dyadic prime and the narrow class group ${A_F^\prime}^+ =0$. 
In particular, the vanishing of  $H^2_+(o_F^\prime, \Z_2(i))$  is independent of the integer  $i\geq 2$. 
\end{proposition}

\textbf{Proof}. 
the  $2$-rank of  $ H^2( o_F^\prime, \Z_2(i))$  is given by  \cite[Proposition 6.13]{Weibel00} or \cite[Lemma 2.2]{Kolster03}: \\
(i) $rk_2 H^2( o_F^\prime, \Z_2(i)) = \mid S_2 \setminus S_\infty \mid + rk_2({A_F^\prime}))-1$  for  $i$  odd;\\
(ii) $rk_2 H^2( o_F^\prime, \Z_2(i)) = \mid S_2 \setminus S_\infty \mid + rk_2({A_F^\prime}))-1+r_1$  for  $i$  even;\\
where  $\mid S_2 \setminus S_\infty \mid $ represents the number of dyadic primes in  $F$ and  $r_1$  is the number of real places in  $F$.

Also, the exact sequence 
$$0 \to H^1( o_F^\prime, \Z_2(i))/2 \to  H^1( o_F^\prime, \Z/2(i)) \to {_2}H^2( o_F^\prime, \Z_2(i)) \to 0$$
shows that   $H^2( o_F^\prime, \Z_2(i)) =0$  precisely  when 
$$H^1( o_F^\prime, \Z_2(i))/2\cong  H^1( o_F^\prime, \Z/2(i)).$$

Besides,  by Kummer theory we have an exact sequence 
%({\bf Trieste page 224,  mais ca devrait etre assez connu et c'est bizarre qu'on n'y a pas pense})
$$0 \to U_F^{\prime}/U_F^{{\prime }^2} \to  H^1( o_F^\prime, \mu_2) \to {_2A_F^{\prime}} \to 0.$$
Therefore,  $A_F^{\prime}$  vanishes precisely when 
 $$H^1( o_F^\prime, \mu_2) \cong  U_F^{\prime }/U_F^{{\prime }^2}.$$
%
%Besides, let  $F_{S_2}$ be the maximal
%extension of  $F$  unramified outside  $2$ and infinity and 
%$G_{S_2}(F)=\Gal(F_{S_2}/F)$. Then by Kummer theory we have an exact sequence 
%$$0 \to U_F^{\prime}/U_F^{{\prime }^2} \to  H^1( G_{S_2}(F), \mu_2) \to {_2H^1( G_{S_2}(F), U^{\prime})} \to 0$$
%where  $U^{\prime} = \varinjlim_{L} U_L^{\prime}$  with  $L$  running over all finite extensions of  $F$  contained in  $F_{S_2}$. Now since  $H^1( G_{S_2}(F), U^{\prime}) \cong c\ell(o^\prime_F)$, 
% the $2$-primary part of the class group $c\ell(o^\prime_F)$ vanishes precisely when 
% $$H^1( G_{S_2}(F), \mu_2) \cong  U_F^{\prime }/U_F^{{\prime }^2}.$$
%
First assume  $i$ is odd. 
The exact sequence (\ref{kolster2}) shows that  $H^2_+(o_F^\prime, \Z_2(i))=0$   precisely when 
$ H^2( o_F^\prime, \Z_2(i)) =0$ and the signature map  (\ref{sgn}) is surjective. 
By the above formula,  $ H^2( o_F^\prime, \Z_2(i))=0 $ if and only if  $F$ contains only one dyadic prime and ${A_F^\prime}=0$,  which by the above discussion is equivalent to  $\mid S_2 \setminus S_\infty \mid=1$  and
$$H^1( o_F^\prime, \Z_2(i))/2\cong  H^1( o_F^\prime, \Z/2(i)) \cong  U_F^{\prime }/U_F^{{\prime }^2}.$$ 
In particular, the  signature map (\ref{sgn}) becomes 
$$\sgn_F : U_F^{\prime }/U_F^{{\prime }^2}  \rightarrow  (\Z/2\Z)^{r_1}$$
whose surjectivity  is equivalent to the equality  ${A_F^\prime}^+ ={A_F^\prime}$ 
(cf. \S 1.2).

 If instead $i$ is even, consider the  following commutative diagram 
$$
\begin{array}{cccccccccccc}
0 & \longrightarrow & H^1(o^\prime _F, {\Z_2}(i))/2& \longrightarrow & H^1(o^\prime _F, \Z/2(i)) & \longrightarrow &
{_2H^2(o^\prime _F, {\Z_2}(i))}& \longrightarrow & 0\\
 &  &  \downarrow  &    & \sgn  \downarrow  &             &  \sigma  \downarrow \\
 &  & 0 & \longrightarrow & ( \Z/2)^{r_1} & \stackrel{\cong} 
{\longrightarrow} &  ( \Z/2)^{r_1}& &  
  \end{array}
  $$
obtained in the same manner as diagram (\ref{diagram sgn})  where the left vertical map is induced by the signature map. 

The short exact sequence (\ref{H2+H2 i even}),  readily yields an exact sequence 
$$0\rightarrow  {_2H^2_+(o_F^{\prime}, \Z_2(i))} \rightarrow {_2H^2(o_F^{\prime},
\Z_2(i))} \stackrel{\sigma} \rightarrow ( \Z/2)^{r_1} \rightarrow H^2_+(o_F^{\prime}, \Z_2(i))/2.$$
Hence, we successively have 

\begin{tabular}{llll}
$H^2_+(o_F^{\prime}, \Z_2(i))=0$ & 
%$\Longleftrightarrow$ &  $ _2H^2(o_F^{\prime}, \Z_2(i)) \stackrel{\sigma} \longrightarrow ( \Z/2)^{r_1}$  is an isomorphism\\
 $\Longleftrightarrow $ & $\sigma$ is surjective and $rk_2 H^2( o_F^\prime, \Z_2(i)) =r_1$\\
& $\Longleftrightarrow $ &   $\mid S_2 \setminus S_\infty \mid=1$, ${A_F^\prime}=0$  and 
$H^1(o^\prime _F, \Z/2(i)) \stackrel{\sgn} \longrightarrow ( \Z/2)^{r_1}$ is surjective \\
 & $\Longleftrightarrow $ &   $\mid S_2 \setminus S_\infty \mid=1$, ${A_F^\prime}=0$  and 
 $ U_F^{\prime }/U_F^{{\prime }^2} \stackrel{\sgn} \longrightarrow ( \Z/2)^{r_1}$ is surjective \\
 & $\Longleftrightarrow $ &   $\mid S_2 \setminus S_\infty \mid=1$  and ${A_F^\prime}^+=0.$ 
\end{tabular}
\qed 

A number field  $F$ which has only one dyadic prime and whose narrow class group ${A_F^\prime}^+ =0$   is called  $2$-regular (\cite{Gras-Jaulent89, Nguyen90,Ostvaer00}). 
According to the above proposition a number field  $F$  is  $2$-regular precisely when the positive cohomology  $H^2_+(o_F^{\prime}, \Z_2(i))$  vanishes for some (hence all) integer  $i$.

Theorem \ref{motivic+ codescent thm}  leads  to the following going-up property which completes Corollarry  \ref{going up p,i regularity}  for  $p=2$:
\begin{cor}\label{going up 2,i regularity}
Let  $L/F$ be  a finite Galois $2$-extension  of number fields  with
Galois  group  $G$. Then  $H^2_+(o_L^\prime, \Z_2(i))=0$  precisely
when  the set of primes of $F$ which are in $S_2$ or ramify in   $L$ is $D_F^{+(i)}$-primitive for $(F,2)$.
\qed \end{cor}

Recall that   $D_F^{+(i)} = D_F^{(i)}$  for  $i$ even according to Proposition  \ref{sgn triviality}. 

%\textbf{Ici la remarque concernant la question de B. Kahn?}

We now focus on the case of cyclic extensions of degree $p$, $p$
prime, and give general genus formulae involving only the number of
some specified tamely ramified primes.
%extending \cite{Assim-Movahhedi04, Assim-Movahhedi12}.
As before let $E=F(\mu_p)$ and  let  $D$  be a subgroup of
$E^\bullet$ such that  $E(\sqrt[p]{D})/E$ is unramified outside
$p$-adic primes. We have a perfect pairing

\begin{equation}
\label{pairing}
\begin{array}{ccc}
{\rm Gal}(E(\sqrt[p]{D})/E) \; \times \; D/{E^{\bullet^p}} &
\longrightarrow & \mu_p \\
(\sigma, a) & \longmapsto &  \sigma (\sqrt[p]{a})/\sqrt[p]{a}.
\end{array}
\end{equation}

Let $T$ be  a maximal $D$-primitive set contained in $S$. For each
$v \in S$,  denote by $\sigma_v$ the Frobenius at the prime $v$ in
$E(\sqrt[p]{D})/E$. Finally, let $H := <\sigma_v,\, v \in S-S_p >=
<\sigma_v,\, v \in T-S_p >$  and let $H^{\perp}$ be the orthogonal
complement of $H$ under the above pairing $(\ref {pairing})$.

We have the following lemma which will be crucial in the proof of Proposition  \ref{Norm index formula}  below.

\begin{lem}
 Under the pairing $(\ref {pairing})$, the orthogonal complement
$H^{\perp}$ of $H$  is equal to  the subgroup
$$D \cap
\bigcap_{v \in S-S_p} N_{L_w/F_v}(L_w^\bullet)/{E^{\bullet ^p}}$$ of
$D/E^{\bullet ^p}$.
\end{lem}

{\bf Proof.} The proof is identical to that of \cite[Theorem
2.4]{Assim-Movahhedi12} (see also \cite[Proposition
3.6]{Assim-Movahhedi04}) where  $D_F^{(i,n)}$  is replaced by  $D$.
\qed

Since
$$[D: D \cap \bigcap_{v \in S-S_p} N_{L_w/F_v}(L_w^\bullet)]= \dim_{F_p} ({\rm Gal}(E(\sqrt[p]{D})/E))- \dim_{F_p} (H^{\perp})$$
and  $\dim_{F_p} (H^{\perp}) =\dim_{F_p} ({\rm
Gal}(E(\sqrt[p]{D})/E))- \dim_{F_p} (H)$, we have the following

\begin{proposition}
\label{Norm index formula} Let $L/F$  be a cyclic extension of
degree  $p$. Then, with the above notation, we have the following
formula for the norm index
$$[D: D \cap
\bigcap_{v \in S-S_p} N_{L_w/F_v}(L_w^\bullet)]=p^{t_D},$$ where
$t_D:=\dim_{\F_p}<\sigma_v(E(\sqrt[p]{D})/E)/ v \in S-S_p>$ denotes the maximal number of  tamely ramified primes in $L/F$
belonging to a $D$-primitive set for ($F, \; p$).\qed
 \end{proposition}

 When $L/F$ is a Galois extension of degree a prime number $p$, the norm index occurring in  the genus formula in Theorem
   \ref{motivic genus} can be written as
$${[D_F^{(i)}: D_F^{(i)} \cap (\bigcap_{v\in S \setminus S_p} N_{L_w/F_v}(L_w^\bullet))]}.$$

Indeed, we have by (\ref{p-adic local residual 1/1}),  an isomorphism
      $$H^1(F_v, {\Zp}(i) ) \cong H^1(k_v, {\Zp}(i) )  $$
for any prime $v\in S \setminus S_p$.
Since $ v\in S \setminus S_p$ is allowed to ramify in a $p$-extension, $F_v$ contains the group $\mu_p$ of $p$-th roots 
of the unity. As in the global case, there exists a subgroup  $D_v^{(i)}$ such that 
$$H^1(F_v, {\Zp}(i))/p \cong D_v^{(i)}/F_v^{\bullet p}.$$
Since $L/F$ is cyclic of degree $p$,  the norm index  in Theorem
   \ref{motivic genus} is the same as the order of the image of the map
$$H^1(F, {\Zp}(i) ) \longrightarrow \oplus_{v\in S \setminus S_p}  H^1(F_v, {\Zp}(i) )/N_{L_w/F_v} (H^1(L_w, {\Zp}(i) )
\cong D_v^{(i)}/ F_v^{\bullet p}N_{L_w/F_v}  D_w^{(i)}. $$
By  \cite[Lemma 2.3]{Assim-Movahhedi12}  (whose proof goes along the same line for $p=2$), 
$$D_v^{(i)}/ F_v^{\bullet p}N_{L_w/F_v}  D_w^{(i)} \cong D_v^{(i)}/ D_v^{(i)}\cap N_{L_w/F_v}  (L_w^\bullet)$$ for any 
$ v\in S \setminus S_p$.

 \begin{cor} \label{genus p} Let  $L/F$  be a cyclic extension of prime degree  $p$.
  Then,
  \begin{enumerate}
 \item If $p$ is odd,
 $$\frac{|(K_{2i-2}o_L)_{G}|}{|K_{2i-2}o_F|}=\frac{|H^2_\M(o_L, {\Z}(i))_{G}|}{|H^2_\M(o_F, {\Z}(i))|}
 = p ^{|S \setminus S_p|-t_i}.$$
 \item If $p=2$ and  $2i-2\equiv 2 \pmod 8$,
 $$\frac{|(K_{2i-2}o_L)_{G}|}{|K_{2i-2}o_F|}=\frac{|H^2_\M(o_L, {\Z}(i))_{G}|}{|H^2_\M(o_F, {\Z}(i))|} =  2^{|S \setminus S_p|-t_i-r}.$$
\item If $p=2$ and $2i-2\equiv 6 \pmod 8$,
 $$\frac{|(K_{2i-2}o_L)_{G}|}{|K_{2i-2}o_F|}= 2^r\frac{|H^2_\M(o_L, {\Z}(i))_{G}|}{|H^2_\M(o_F, {\Z}(i))|} =  2^{|S \setminus S_p|-t_i}.$$
 \end{enumerate}
 Here  $t_i :=\dim_{\F_p}<\sigma_v(E(\sqrt[p]{D_F^{(i)}})/E)/ v \in S-S_p>$ denotes the maximal number of  tamely ramified primes in
$L/F$ belonging to a $D_F^{(i)}$-primitive set for ($F, \; p$)  and
$r$  is the number of infinite places of  $F$ which ramify in  $L$.
 \end{cor}

{\bf Proof}. By Proposition \ref{Norm index formula}, in all the
considered cases, we have
$${[D_F^{(i)}: D_F^{(i)} \cap (\bigcap_{v\in S \setminus S_p} N_{L_w/F_v}(L_w^\bullet))]} =p^{t_i}.$$
The formulae for the motivic cohomology groups follow from Theorem
\ref{motivic genus}. When  $p$ is odd or  $p=2$  and   $2i-2\equiv 2
\pmod 8$, the formulae for the $K$-groups is then a consequence of
Theorem 3.1. In the remaining case, we use the formula (\ref{K6
motivic}). \qed

When  $p=2$ and $i\geq 2$  is odd, we use positive cohomology. The norm
index in Theorem   \ref{genus+} can be written as
$${[D_F^{+(i)}:D_F^{+(i)} \cap (\bigcap_{v\in S \setminus S_p} N_{L_w/F_v}(L_w^\bullet))]}$$
and recall that we have an exact sequence
$$0\to  D_F^{+(i)}/F^ {\bullet 2} \to  D_F^{(i)}/F^ {\bullet 2}  \to (\Z/2)^{r_1} \to (\Z/2)^{\delta_i} \to 0$$
where the second map is the restriction of the partial signature map
$\sgn_{R_{L/F}}$  to  $D_F^{(i)}/F^ {\bullet 2}$.

%{\bf Peut-on comprendre la valeur de  $s_i$  dans ce cas quadratique  $L/F$}

\begin{cor} \label{genus 2 i odd}
Let  $L/F$  be a quadratic extension of number fields and $i\geq 2$ odd. 
Let  $t_i^+$ denote the maximal number of  tamely
ramified primes in $L/F$ belonging to a $D_F^{^+(i)}$-primitive set
for ($F, \; 2$).  Then
$$ \frac{|H^2_+(o'_L, {\Z}_2(i))_{G}|}{|H^2_+(o'_F, {\Z}_2(i))|}
= 2 ^{|S \setminus S_2|-t_i^+}.$$ 
Moreover if hypothesis $(\mathcal{H}_i)$ holds, then 
$$\frac{|H^2(o'_L, {\Z_2}(i))_{G}|}{|H^2(o'_F, {\Z_2}(i))|} =  2 ^{|S \setminus S_2|+s_i-t_i^+}$$ 
and
\begin{enumerate}
\item If  $2i-2\equiv 0 \pmod 8$,
 $$\frac{|(K_{2i-2}o_L)_{G}|}{|K_{2i-2}o_F|}=\frac{|H^2_\M(o_L, {\Z}(i))_{G}|}{|H^2_\M(o_F, {\Z}(i))|} =  2^{|S \setminus S_2|+s_i-t_i^+}.$$
\item If  $2i-2\equiv 4 \pmod 8$,
 $$\frac{|(K_{2i-2}o_L)_{G}|}{|K_{2i-2}o_F|}= 2^r\frac{|H^2_\M(o_L, {\Z}(i))_{G}|}{|H^2_\M(o_F, {\Z}(i))|} =  2^{|S \setminus S_2|-t_i^+}.$$
 \end{enumerate}
 \end{cor}
 
{\bf Proof}. By Proposition \ref{Norm index formula}, we have
$${[D_F^{+(i)}: D_F^{+(i)} \cap (\bigcap_{v\in S \setminus S_p} N_{L_w/F_v}(L_w^\bullet))]} =p^{t_i^+}.$$
By  Theorem \ref{genus+}, it follows that
$$\frac{|H^2_+(o'_L, {\Z}_2(i))_{G}|}{|H^2_+(o'_F, {\Z}_2(i))|} = 2 ^{|S \setminus S_2|-t_i^+}.$$
The formulae for the motivic cohomology groups then follow from
Proposition \ref{compare + motivic}. Now, if  $2i-2\equiv 0 \pmod
8$,  then the formula follows by  Theorem  \ref{RW} and if
$2i-2\equiv 4 \pmod 8$, then we use formula (\ref{comparing K4 genus
motivic genus}). \qed

When $2i-2\equiv 0 \pmod 8$ and  $L/F$ is unramified at infinite primes,  we use  
Theorem \ref{odd cyclic genus} to obtain an explicit formula without any assumption on the cohomology of
$(Z/2)^{\delta_i(L)}$:

\begin{proposition} \label{genus 2 i odd ur}
Let  $L/F$  be a quadratic extension of number fields and $i\geq 2$ odd. 
Let  $t_i$ denotes the maximal number of  tamely
ramified primes in $L/F$ belonging to a $D_F^{(i)}$-primitive set
for ($F, \; 2$). Assume that $L/F$ is unramified at infinity. 
Then
$$\frac{|H^2_\M(o_L, {\Z}(i))_{G}|}{|H^2_\M(o_F, {\Z}(i))|} =  2^{|S \setminus S_2|-t_i}.$$
In particular if  $2i-2\equiv 0 \pmod 8$,
$$\frac{|(K_{2i-2}o_L)_{G}|}{|K_{2i-2}o_F|}= 2^{|S \setminus S_2|-t_i}.$$
\qed 
\end{proposition}

We now have the ingredients to answer a question (of independent interest), raised by B. Kahn in \cite[page 2]{Kahn97}:  
Let $\rho_i:= \rho_i(F)$ denote  the  $2$-rank of  the image of the signature map
$$\mathrm{sgn}_F : H^1(F, {\Z_2}(i))/2 \longrightarrow  ({\Z}/2)^{r_1}.$$ 
Is it true that $\rho_i =1$, for any number field $F$ such that $r_1\geq 1$?
The answer turns out to be negative in general. For instance,  $\Q(\sqrt{5}) $ is $2$-regular by Proposition \ref{description 2,i regularity}, and we saw in the proof that the signature map for such a field is surjective. 
In fact, as the following theorem shows, the integer  $\rho_i $  can be arbitrary large for real number fields. 

\begin{thm} \label{kahn question} Let $n\geq 1$ be an integer.  Then there exists a totally real number field $F$ such that the image of  
the signature map
$\mathrm{sgn}_F$ has  $2$-rank $\rho_i=n$, for all $i\geq 2$ odd.
\end{thm}
 \textbf{Proof.}
Let  $m$  be the unique  integer such that
$$2^m \leq n <2^{m+1}.$$
Since $\Q$ is $2$-regular, there  exist  $2$-regular number fields $F$ such that $[F: \Q]= 2^m$  by Theorem \ref{motivic+ codescent thm}. 
For instance, one may take for  $F$  the  $m$-th layer of the cyclotomic  $\Z_2$-extension of  $\Q$. 

Then, as noticed in the proof of Proposition \ref{description 2,i regularity}, the signature map  $\sgn_F$  is surjective and the exact sequence   (\ref{sequence D+ and D}) shows that 
$$ \dim_{\mathbf{F}_2}   D_F^{+(i)}/F^ {\bullet 2} = \dim_{\mathbf{F}_2}  D_F^{(i)}/F^ {\bullet 2}  -  r_1(F).$$
% The exact sequence   (\ref{sequence D+ and D}) shows that  for any number field $F$, 
%$$ \dim_{\mathbf{F}_2}   D_F^{+(i)}/F^ {\bullet 2} = \dim_{\mathbf{F}_2}  D_F^{(i)}/F^ {\bullet 2}  -  \rho_i.$$
%Therefore, if $\mathrm{sgn}_F$ is surjective, then 
%$$ \dim_{\mathbf{F}_2} D_F^{(i)}/F^ {\bullet 2} =  \dim_{\mathbf{F}_2}   D_F^{+(i)}/F^ {\bullet 2}  - r_1.$$
%This is the case if $F$ is $2$-regular, as seen in the proof of Proposition  \ref{description 2,i regularity}.
%Assume now that  $F$ is a  $2$-regular totally real number field.  In particular, $\rho_i(F)= r_1$.
Choose a maximal $D_F^{+(i)}$-primitive set  $T$  for ($F,2$)   (Definition \ref{defi primitive}). 
The set  $T$  then contains the  dyadic prime and exactly one non dyadic prime but we do not need this fact.     

Let $S$ be a  $D_F^{(i)}$-primitive set for ($F,2$)  containing $T$  and such that   $|S\setminus T | = 2^{m+1}-n$. 
Such a primitive set  $S$  exists since  $2^{m+1}-n \leq 2^{m+1}-2^m =2^m =  r_1(F).$ 
Then there exists a real quadratic extension  $L$ of $F$  such that the non dyadic primes which are ramified in  $L/F$ consist precisely of 
$S\setminus S_2$   (see lemma \ref{existence} below).
Proposition \ref{genus 2 i odd ur} shows that 
$H^2(o_L^\prime, {\Z_2}(i))=0$. By the exact sequence (\ref{s_i as kernel}), we then have  
$$H^2_+(o_L^\prime, {\Z_2}(i))\cong  (\Z/2)^{\delta_i(L)}.$$
 On the other hand, Corollary \ref{genus 2 i odd} shows that  
 $$|H^2_+(o_L^\prime, {\Z_2}(i))_G|=2^{|S\setminus T |}$$ leading to   
$\dim_{\mathbf{F}_2} H^2_+(o_L^\prime, {\Z_2}(i)) = |S\setminus T |.$  
Hence  $\delta_i(L) = |S \setminus T | = 2^{m+1}-n$  and   
$$\rho_i(L)= r_1(L)- \delta_i(L) =  2^{m+1}- 2^{m+1}+n = n.$$
%Since  $S$  can be chosen so that   $|S\setminus T | $  takes all values between $0$ and $r_1$,  
%the same holds  for $\delta_i(L)$.  It follows that 
%we can choose $\rho_i(L)= r_1(L)- \delta_i(L) =2r_1- \delta_i(L)$   to be any integer between  $ r_1$  and  $2r_1$:
%$$ r_1 \leq \rho_i(L)  \leq 2r_1.$$
%Now let $n\geq 1$, there exists a unique  integer $m$ such that
%$$2^m \leq n <2^{m+1}.$$
%Since $\Q$ is $2$-regular, the existence of  $2$ regular number field $F$ such that $[F: \Q]= 2^m$ is guaranteed by Theorem \ref{motivic+ codescent thm}. 
%For instance, one may take the  $m$-th layer of the cyclotomic  $\Z_2$-extension of  $\Q$. 
It remains to prove:

\begin{lem} \label{existence}  Let $F$ be a number field with a trivial $2$-primary narrow class group and let $\{\mathfrak{p}_1, \cdots ,\mathfrak{p} _t\}$ be a set 
of non dyadic primes.  Then there exists a (totally real) quadratic extension $L$ of $F$  unramified at infinity and 
in which the tamely ramified primes are precisely  $\mathfrak{p}_1, \cdots ,\mathfrak{p} _t$. 
%which is ramified at $\mathfrak{p}_1, \cdots ,\mathfrak{p} _t$ and unramified outside
%$\{\mathfrak{p}_1, \cdots ,\mathfrak{p} _t\}\cup S_2.$ 
\end{lem}

\textbf{Proof.}  Since the $2$-primary part of the narrow class group of $F$ is trivial, there exists an odd integer  $m$  such that $\mathfrak{p}_i^m = (d_i)$ is a principal ideal generated by a totally positive element for all $i$, $1\leq i\leq t$. 
Obviously  $d_i \not \in F^2$ and  
%one proves by Kummer theory  that 
the field  $L= F(\sqrt{d_1 \cdots d _t})$ 
%is  ramified at $\mathfrak{p}_1, \cdots ,\mathfrak{p} _t$  and unramified outside  
%$\{\mathfrak{p}_1, \cdots ,\mathfrak{p} _t\}\cup S_2.$ 
fulfils the conditions of the lemma. 
\qed

\section{Galois descent}
We keep the notations of the preceding sections: $L/F$ is a finite
Galois extension of number fields with Galois group  $G$. Recall
that  $S$  consists of finite primes which ramify in  $L/F$ as well
as the infinite primes. Suppose that  $S$  contains a  set  $T$ such
that for all prime  $p$ (dividing $e$) the set $T_{L/F,p}:= \{v \in
T /\;  p\mid e_v\}\cup S_p$ is $D_F^{(i)}$-primitive for $(F,p)$.

As before, the natural map
\begin{equation}\label{betaT}
\beta_T : \hat{H}^0(G, H^1_\M(L,\Z(i))) \rightarrow \hat{H}^0(G,
\oplus_{v\in T, w\mid v} H^1_\M(k_w,\Z(i)))
\end{equation}
is surjective.

Recall that, by Theorem \ref{galois descent and codescent},  the
kernel and cokernel of the functorial map
$$  f_i : H^2_\M(F, {\Z}(i))  \longrightarrow H^2_\M(L, {\Z}(i))^G $$
are described by the  $G$-cohomology of  $H^1_\M(L,\Z(i))$.

We shall use the surjectivity of  $\beta_T$  to prove that of the
natural map
$$ \theta_T : {H}^2(G, H^1_\M(L,\Z(i))) �\rightarrow {H}^2(G, \oplus_{v\in T, w\mid v} H^1_\M(k_w,\Z(i))) $$
by taking the cup-product with  $H^2(G, \Z)$. This provides lower
bounds of respectively the order and the generator rank of
${\mathrm{coker}} f_i.$
%Note that by Borel's theorem on the $\Z$-structure  odd  $K$-groups, it is not hard to give upper bounds for the order or the generator rank of
%${\mathrm{coker}} f_i$.
%For each prime $v$ in $T$, fix a prime $w$ of  $L$ above $v$. Let $G_v$, $I_v$ and  $I_v^\prime$ be as in Lemma 2.2.
%Introduce $G_{v}^{\prime}:= G_v/I_v ^{(\ell)}$.

Since $\ell$ (the residue characteristic of $F_v$) does not divide
the order of $H^1_\M(k_w,\Z(i))$, we have
$$ \hat{H}^q(I_v, H^1_\M(k_w,\Z(i))) \cong  \hat{H}^q(I_v^\prime, H^1_\M(k_w,\Z(i)))$$
for all $q\in \Z$.  Here  $I_v^\prime = I_v/I_{v,1}$  is the inertia
group of  $v$  in $L/F$  factored by the first ramification group.

Since $I_v^\prime$ is a cyclic group, cup product gives rise to the
following commutative diagram
$$
  \begin{array}{cccccc}
    \hat{H}^0(I_v^\prime, H^1_\M(k_w,\Z(i)))   &  \otimes    &   {H}^2(I_v^\prime, \Z)   &  \stackrel{\sim}  \rightarrow  &
     {H}^2(I_v^\prime, H^1_\M(k_w,\Z(i)))\\

    \mathrm{cor} \;  {\downarrow}\wr                                            &                  & \mathrm{ res} \; \uparrow &
     & \mathrm{cor} \;  {\downarrow}\wr \\
      \hat{H}^0(I_v, H^1_\M(k_w,\Z(i)))                     & \otimes &   {H}^2(I_v, \Z) & \rightarrow  &   {H}^2(I_v, H^1_\M(k_w,\Z(i)))\\
  \end{array}
  $$
which shows the surjectivity of the bottom map. In fact this map is
an isomorphism since $ H^1_\M(k_w,\Z(i))$  is of order prime to
$\ell$  and the restriction map ${H}^2(I_v, \Z) \to
{H}^2(I_v^\prime, \Z)$  is an isomorphism on the non-$\ell$-part.

Now consider the commutative diagram
$$
\begin{array}{cccccc}
\hat{H}^0(I_v, H^1_\M(k_w,\Z(i))) &  \otimes & {H}^2(I_v, \Z) &   \twoheadrightarrow  & {H}^2(I_v, H^1_\M(k_w,\Z(i)))\\
\mathrm{cor} \;  {\downarrow} & &  & & \mathrm{cor} \;  {\downarrow} \\
\hat{H}^0(G_v, H^1_\M(k_w,\Z(i))) &  & \mathrm{ res} \; \uparrow & & {H}^2(G_v, H^1_\M(k_w,\Z(i)))\\
{\downarrow} \wr & & & &   {\downarrow} \wr\\
\hat{H}^0(G, \oplus _{w\mid v}H^1_\M(k_w,\Z(i))) & \otimes & {H}^2(G, \Z) & \stackrel{\gamma_v} \rightarrow  & {H}^2(G, \oplus _{w\mid v}  H^1_\M(k_w,\Z(i))). \\
  \end{array}
  $$

 By (\ref{residual shifting}) the cokernel of the map
  \begin{equation}\label{cor H2 Iv Gv (i)}
\mathrm{cor}:  {H}^2(I_v, H^1_\M(k_w,\Z(i)))  \rightarrow {H}^2(G_v,
H^1_\M(k_w,\Z(i)))
  \end{equation}
   is isomorphic to the cokernel of the corestriction
  \begin{equation}\label{cor H0 Iv Gv (i-1)}
  \hat{H}^0(I_v, H^1_\M(k_w,\Z(i-1)))  \rightarrow   \hat{H}^0(G_v, H^1_\M(k_w,\Z(i-1)))
    \end{equation}
induced by the norm. The map (\ref{cor H0 Iv Gv (i-1)})  is surjective
since, as in the proof of Lemma \ref{lemma evi}, we have
$\hat{H}^0(I_v, H^1_\M(k_w,\Z(i-1))) \cong H^1_\M(k_w,\Z(i-1))/e_v$
and $\hat{H}^0(G_v, H^1_\M(k_w,\Z(i-1))) \cong
H^1_\M(k_v,\Z(i-1))/e_v.$ Hence the map (\ref{cor H2 Iv Gv (i)})  is also
surjective. It follows from the above diagram that
    the  cup-product induces a surjective homomorphism
  $$ \gamma_v :       \hat{H}^0(G, \oplus _{w\mid v}H^1_\M(k_w,\Z(i)))   \otimes
   {H}^2(G, \Z)  \rightarrow    {H}^2(G, \oplus _{w\mid v} H^1_\M(k_w,\Z(i))).
$$
Consider now the following commutative diagram
  $$
  \begin{array}{ccccccccc}
\hat{H}^0(G, H^1_\M(L,\Z(i)))   \otimes {H}^2(G, \Z)       &
\rightarrow      &
{H}^2(G, H^1_\M(L,\Z(i))) \\
     {\downarrow}                                                                                      &                        &  \downarrow \theta_T \\
\hat{H}^0(G,  \oplus_{v\in T, w\mid v} H^1_\M(k_w,\Z(i)))   \otimes
{H}^2(G, \Z) &     \stackrel{\oplus \gamma_v} \twoheadrightarrow  &
       {H}^2(G, \oplus_{v\in T, w\mid v} H^1_\M(k_w,\Z(i)))  \\
  \end{array}
  $$
where the horizontal  maps are induced  by cup-product.  The left
vertical map is induced by $\beta_T$ (see (\ref{betaT})) and so is
surjective.  The surjectivity of the right vertical map  $\theta_T$
% $$ {H}^2(G, H^1_\M(L,\Z(i)))  \rightarrow {H}^2(G, \oplus_{v\in T, w\mid v} H^1_\M(k_w,\Z(i))) $$
 then follows from the surjectivity of the maps $\gamma_v$.
 To compute the cohomology group\\
$H^2(G, \oplus_{v\in T, w\mid v} H^1_\M(k_w,\Z(i)))$, we have for
$v\in T$,
 $$
 \begin{array}{ccl}
   H^2(G_v,  H^1_\M(k_w,\Z(i))) & \cong &   \hat{H}^0(G_v,  H^1_\M(k_w,\Z(i-1))) \;
   \quad (\mbox{by (\ref{residual shifting})})\\
    & \simeq &  \Z/e_{v}^{(i-1)}  \; \quad (\mbox{see the proof of Lemma \ref{lemma evi}})
 \end{array}
 $$
 where $e_{v}^{(i-1)} =gcd(e_v, q_v^{i-1}-1)$.

 The above discussion together with Theorem \ref{galois descent and codescent} yield the following upper bounds for the capitulation cokernel:
 \begin{thm}\label{bounds coker}
Let  $L/F$ be a Galois  extension  of number fields   with Galois
group  $G$. Assume that the set of ramified primes in the extension
$L/F$ contains a set $T$ such that for all prime $p$, the set
$\{v\in T/ \; p\mid e_v\} \cup S_p$ is $D_F^{(i)}$-primitive for
$(F,p)$. Then we have a surjective homomorphism
 $$H^2(G, H^1_\M(L,\Z(i)))  \longrightarrow  \oplus_{v\in T}  \Z/e_{v}^{(i-1)}$$
leading to the following upper bounds:
\begin{enumerate}
\item $\mid {\mathrm{coker}} f_i\mid \geq \prod_{v \in T}e_{v}^{(i-1)}$  for  $i$  even and \\
\item $\mid {\mathrm{coker}} f_i\mid \geq 2^{s_i-r}\prod_{v \in T}e_{v}^{(i-1)}$  for  $i$  odd
\end{enumerate}
 where for each $v$, $e_{v}^{(i-1)} :=gcd (e_v, q_v^{i-1}-1)$.
 \qed
\end{thm}
When $G$ is a cyclic group, we see from Theorem \ref{galois descent
and codescent}  that the Herbrand quotient $h(G,  H^1_\M(L,\Z(i)))$
equals $2^{-r}$ for  $i$ even and $2^r$ for  $i$ odd. Hence the
above theorem together with Lemma \ref{lemma evi} lead to the
following upper bounds for the capitulation kernel:
 \begin{cor}\label{bounds ker}
Let  $L/F$ be a cyclic extension  of number fields of degree  $ n$
with Galois  group  $G$.
 Assume that the set of ramified primes in the extension $L/F$ contains
a set $T$ such that for all prime $p$, the set $\{v\in T/ \; p\mid
e_v\} \cup S_p$ is $D_F^{(i)}$-primitive for $(F,p)$. Then

 $\mid\ker f_i\mid   \geq 2^r \prod_{v \in T}e_v^\prime$ for $i$ even and
 $\mid\ker f_i\mid   \geq 2^{-r}\prod_{v\in T}e_v^\prime$ for $i$ odd.
\qed
\end{cor}
We finish the section by giving examples of Galois extensions $L/F$
for which the kernel and the cokernel of  $f_i$ are explicitly
known.
%We say that a number field $F$ is $(n,i)$-regular if the  group $H^2(o_F,
%\Z(i))/n$ vanishes, or equivalently $F$ is $(p,i)$-regular for all
%prime $p$ dividing $n$. For example, one of the equivalent
%formulations of Vandiver's conjecture is that the field $\Q$ of
%rational numbers is  $(n,i)$-regular for odd $i$ and all integers
%$n$ \cite{Kurihara92}. Moreover, for any number field $F$, the finiteness of
%  $H^2(o_F, \Z(i))$ shows that $F$ is  $(n,i)$-regular for an infinite number of integers    $n$.

Since $H^2_\M(o_F, \Z(i))$ is finite there exist infinitely many
integers $n$  such that  $H^2_\M(o_F, \Z(i))/n$ vanishes. We take
such an integer  $n$  and a cyclic extension  $L/F$ of degree $n$,
with group $G$ and unramified at archimedean places. As before, $S$
consists of ramified primes as well as archimedean ones. Assume that
for all prime numbers  $p$, the set $\{v\in S/ \; p\mid e_v\} \cup
S_p$ is $D_F^{(i)}$-primitive for  $(F,p)$. Consider   the
commutative diagram obtained from the exact sequence  (\ref{tame
kernel}):
$$
 \begin{array}{ccccccccccc}
 0 &  \rightarrow & H^2_\M(o_L, \Z(i))^G  & \rightarrow & H^2_\M(o_L^S, \Z(i))^G  &
 \rightarrow & \oplus_{v\in S\setminus S_\infty }({\oplus_{w \mid v}}  H^1_\M(k_w, {\Z}(i-1) ))^G  \\
 & & f_i^\prime  \uparrow & & f_i \uparrow &   &  \oplus f_{v,i}  \, \; \uparrow\\
  0 &   {\rightarrow} &       H^2_\M(o_F, \Z(i))  & {\rightarrow} &  H^2_\M(o_F^S, \Z(i))  & {\rightarrow} &
\oplus_{v\in S\setminus S_\infty }  H^1_\M(k_v, {\Z}(i-1))&
{\rightarrow} & 0
 \end{array}
 $$
By Theorem \ref{galois descent and codescent}, $\ker f_i \simeq
H^1(G, H^1_\M(L,\Z(i)))$  is annihilated by  $n$. Hence the
vanishing of  $ H^2(o_F, \Z(i))/n$ shows that  $f_i^\prime$  is
injective. Moreover, by Theorem \ref{motivic codescent thm},
$H^2_\M(o_L, \Z(i))^G$  has the same order as $H^2_\M(o_F, \Z(i))$.
Therefore $f_i^\prime$  is an isomorphism.

Now,
$$\coker f_{v,i} \simeq \hat{H}^0(G_v, H^1_\M(k_w,\Z(i-1))) \simeq \Z/e_{v}^\prime \quad  \mbox{(see the proof of Lemma \ref{lemma evi})}$$
and similarly
$$\ker f_{v,i} \simeq \hat{H}^{-1}(G_v, H^1_\M(k_w,\Z(i-1))) \simeq \Z/e_{v}^\prime .$$
This leads to the exact structure of the kernel and the cokernel of
$f_i$:
$$ker f_i \simeq {\mathrm{coker}} f_i  \simeq \oplus_{v\in S\setminus S_\infty} \Z/e_{v}^\prime .$$

%\begin{rem}
%As noticed before, the deviation between motivic cohomology and
%$K$-theory is known  \cite{Kahn97, Weibel00}. When $2i-2\equiv 0, 2
%\pmod 8$, it is known that the motivic Chern characters
%$$
%    ch_{i,2}^\M : K_{2i-2} (F)   \longrightarrow  H^2_\M (F, {\Z}(i)).
%   $$
%are isomorphisms.  In that cases, the above results are also valid
%for the higher $K$-groups. In the other cases  ($2i-2 \equiv 4, 6
%\pmod 8$), the above statements remain valid for the higher
%$K$-groups up to a power of  $2$.
%\end{rem}
%Let us finish this section with a few comments on the   surjectivity of the canonical map
%$$H^1(F,\Z(i))/n  \longrightarrow   \oplus_{v\in S} H^1(F_v,\Z(i))^0/n$$
%which is the same is the surjectivity of the maps
%
%$$H^1(F,\Z_p(i))/p  \longrightarrow   \oplus_{v\in S \setminus S_p } H^1(F_v,\Z_p(i))/p$$
%for all prime dividing $p$.
%
%Denote by $d_i$ the $\Z$-rank of $H^1(F,\Z(i))$ : $d_i=...$ by Borel's Results.
%
%
%
%Let $E:= F(\mu_p$.  Then there exists a subgroup $D_F^{(i,p)}$ of $E$ such that
%$$H^1(F,\Z(i))/p \cong H^1(F,\Z_p(i))/p \cong D_F^{(i,p)}/ E^{\bullet , ^}$$

\section{examples}
In this section, we consider the special case where $F= \Q$. Given a
prime  $p$ we are going to caracterize  $p$-extensions  $L$  of $\Q$
for which $H^2_\M(o_L, \Z(i))\otimes \Z_p \simeq H^2_\et(o_L^\prime,
\Z_p(i)) =0$. Recall that the motivic cohomology groups $H^2_\M(o_L,
\Z(i))$  and  the K-groups  $K_{2i-2} o_L$  are isomorphic up to
known $2$-torsion  \cite{Kahn97, Weibel00}. Hence for  $p$  odd, we
actually determine  $p$-extensions   $L$  of  $\Q$ for which
${K}_{2i-2}o_L \otimes \Z_p$  vanishes. The list of $p$-extensions
$L$  of  $\Q$  for which  the $p$-part of the classical wild kernel
$WK_2L\otimes \Z_p$  vanishes  can be found in \cite[example 2.13,
2]{Kolster-Movahhedi00}

%For a detailed account on the $K$-theory of rings of integers see for instance \cite{Weibel05}.
\subsection  { }
We begin with the case where  $p$  is an  {\it odd}  prime. In this
case  $H^2_\et(o_L^\prime, \Z_p(i))$  is isomorphic to
${K}_{2i-2}o_L \otimes \Z_p$.

First, consider the case where  $i =2k$  is even. The  order of
$H^2_\M(\Z, \Z(i))$ can be computed using the values of the Riemann
zeta function at odd negative integers since the Lichtenbaum
conjecture is true in this case thanks to the main theorem in
Iwasawa theory (Theorem of Mazur-Wiles) \cite{MW84}. Namely, let
$B_k$  be the  $k$-th Bernoulli number and  $c_k$  be the numerator
of  $B_k/4k$. Then  $\mid H^2_\M(\Z, \Z(i)) \mid = 2c_k$.
For instance, for  $p$ odd and  $i=2, 4, 6, 8,  10$, we have $H^2_\et(\Z[1/p], \Z_p(i))= 0$    whereas for  $p=691$  and  $i=12$, the group $H^2_\et(\Z[1/p], \Z_p(i))$ is cyclic of ordrer $p$.\\

%First, consider the case where  $i =2k$  is even. \\
%The  order of  $K_{2i-2}\Z$
%can be computed using the values of the Riemann zeta function at odd negative integers since the Lichtenbaum conjecture is true in this case thanks to the main theorem in Iwasawa theory (Theorem of Mazur-Wiles).
%Namely, let  $B_k$  be the Bernoulli number and  $c_k$  be the numerator of  $B_k/4k$.
%Then  $\mid K_{2i-2}\Z \mid = c_k$  for  $k$  even and  $= 2c_k$  for  $k$  odd.
%In particular for  $p$ odd and  $i=2, 4, 6, 8,  10$, we have ${K}_{2i-2}\Z \otimes \Z_p= 0$    whereas for  $p=691$  and  $i=12$, we have ${K}_{22}\Z$ is cyclic of ordrer $p$.\\

If   $i \not\equiv 0  \pmod {p-1}$, Corollary \ref{going up p,i
regularity}  shows that $H^2_\et(o_L^\prime, \Z_p(i))=0$ for a
$p$-extension $L$ of $\Q$, precisely when  $H^2_\et(\Z[1/p],
\Z_p(i))= 0$ and  $L$
 is contained in the cyclotomic
$\Z_p$-extension of $\Q$.

If $i \equiv 0  \pmod {p-1}$, the triviality of $H^2_\et(o_L^\prime,
\Z_p(i))$ is equivalent to the  $p$-rationality of the field $L$ and
it is known that a $p$-extension $L$ of $\Q$ is $p$-rational exactly
when at most one non-$p$-adic prime  $\ell$ is ramified in $L$
where  $\ell$ is  such that $\ell \equiv 1 \pmod p$ and  $\ell
\not\equiv 1 \pmod {p^2}$ \cite{Movahhedi88, Movahhedi90,
Movahhedi-Nguyen90}. The last assertion also follows readily from
Corollary \ref{going up p,i regularity}.

%The structure of
%$K$-groups of $\Z$ depends on the validity of  Vandiver's
%conjecture which asserts that the $p$-primary part  of the
%class-group of the maximal real sub-field $\Q(\mu_p)^+$ de
%$\Q(\mu_p)$ is trivial (see \cite{Weibel05} for a detailed account
%on the $K$-theory of $\Z$).\\
%
%By a well-known
%results of Soul\'{e} (up $2$ and $3$-torsion) and Rognes,
%${K}_{4}\Z=0$. It is also known that
% for $p>>0$ and $i$ odd \cite{Soule99},
%${K}_{2i-2}\Z \otimes \Z_p=0$.

Now, consider the case where  $i =2k+1$  is odd. \\
Assuming the triviality of   $H^2_\et(\Z[1/p], \Z_p(i))$ (which is
the case under Vandiver's conjecture), we shall construct all
$p$-extensions $L$
of $\Q$ such that  $H^2_\et(o_L^\prime, \Z_p(i))=0$. \\
Let $E=\Q(\mu_p)$   and  $\Delta=\Gal (E/\Q)$. We have
$\mathrm{dim}_{\mathbb{F}_p}H^1_\et(\Z[1/p], \Z_p(i))/p=1$  and
$$H^1_\et(\Z[1/p], \Z_p(i))/p  \simeq U^\prime_E/p(i-1)^\Delta$$
where  $U^\prime_E$  denotes the group of units of $o_E^\prime$.

If $i\equiv 1 \pmod {p-1}$, $H^1_\et(\Z[1/p], \Z_p(i))/p \simeq
U^\prime_\Q/p(i-1)$ and therefore  $D_Q^{(i)}/E^{\bullet p}$  is
generated by the class of $p$. Hence  $D_Q^{(i)}$-primitive sets for
$(\Q, p)$  are of the form $S_p$ or $T=S_p \cup \{\ell\}$ with
$\ell$  inert in  $E(\sqrt[p]{p})/E$. Therefore $H^2_\et(o_L^\prime,
\Z_p(i)) =0$  when $L/\Q$  is unramified outside  $T=S_p \cup
\{\ell\}$  with $\ell \equiv 1 \pmod p$ and  $p \not\in (\Z/\ell)
^{\bullet p}$.

If  $i\not \equiv 1 \pmod {p-1}$, the hypothesis $H^2_\et(\Z[1/p],
\Z_p(i))=0$ shows that $c\ell_E/p(i-1)^\Delta$  is trivial where
$c\ell_E$  is the class group of  $E$. Hence
$$H^1_\et(\Z[1/p], \Z_p(i))/p \cong U^\prime_E/p(i-1)^\Delta = U_E/p(i-1)^\Delta \cong
C_E/p(i-1)^\Delta ,$$ where  $C_E$ is the group of cyclotomic units
(for the last isomorphism see e.g. \cite[chapter 15, Theorem
15.7]{Washington97}). Let   $\omega$ be the Teichm\"{u}ller
character. Then, for the even integer $j:=1-i$, the class of the
cyclotomic element
$$\xi_j:=\prod_{\delta \in \Delta} (\zeta^\delta - 1)^{\omega^{-j}(\delta)}$$
generates  $D_\Q^{(i)}/E^{\bullet p}$. Hence  $D_\Q^{(i)}$-primitive
sets for  $(\Q, p)$  are of the form $S_p$ or $T=S_p \cup \{\ell\}$
with  $\ell$  inert in $E(\sqrt[p]{\xi_j})/E$. Therefore
$H^2_\et(o_L^\prime, \Z_p(i)) =0$  when $L/\Q$  is unramified
outside  $T=S_p \cup \{\ell\}$  with $\ell \equiv 1 \pmod p$ and
$\xi_j\not\in \Q_\ell^{\bullet p}$.

\subsection  { }
 Assume now $p=2$ and consider a $2$-extension  $L$ of $\Q$. Let,
as before, $S$  be the set consisting of infinite primes as well as
those which  ramify in $L$. The triviality of the class group of
$\Q$  together with Kummer theory show that
\begin{equation} \label{TatekernelforQ}
H^1_\et (\Z[1/2],  {\Z}/2(i)) \cong  U^\prime_\Q/U^{\prime 2}_\Q
=<-{\overline 1},{\overline 2}>
\end{equation}
where, as before,   $U^\prime_\Q$  denotes the group of units in
$\Z[1/2]$.

(i)  If $i$ is even,  Proposition \ref{beginning genus} shows that
$H^2_\et(o_L^\prime,  \Z_2(i))=0$ precisely  when $L$ is totally
imaginary and $\beta_S$ is surjective.  By Lemma \ref{betas loc},
$\beta_S$ is surjective when  $S\cup S_2$ is $D_\Q^{(i)}$-primitive
  for $(\Q, 2)$. Recall that  $S_2=\{ 2, \infty\}$.
  In this case  the signature map is trivial (see \cite[Lemma 2.5]{Kolster03} or Proposition  \ref{sgn triviality} above) so that
 $D_\Q^{(i)}/\Q^{\bullet 2}$  is generated  by the class of $2$.
It follows that $D_\Q^{(i)}$-primitive sets for $(\Q, 2)$  are
precisely those contained in $\{2, \infty, \ell \}$ with $\ell $
inert in  $\Q(\sqrt{2})/\Q$, or equivalently $\ell \equiv \pm 3
\pmod 8$. Therefore
%, when  $i$ is even,  $H^2_\et(o_L^\prime,
%\Z_2(i))=0$  for a  $2$-extension  $L$  of  $\Q$  precisely when $L$
%is totally imaginary and there exists an odd prime $\ell \equiv \pm
%3 \pmod 8$  such that  $L/\Q$  is unramified outside  $\{2, \infty,
%\ell \}$.

\begin{proposition}
Let   $L$ be a finite Galois $2$-extension  of $\Q$ and $i\geq 2$ even. Then the \'{e}tale cohomology $H^2_{\mbox{\'{e}t}}( o_L[1/2],
\Z_2(i))$ vanishes 
precisely when $L$
is totally imaginary and there exists an odd prime $\ell \equiv \pm
3 \pmod 8$  such that  $L/\Q$  is unramified outside  $\{2, \infty,
\ell \}$ with
$\ell \equiv \pm 3 \pmod 8$. \qed
\end{proposition}

%\begin{cor}
%Let   $L$ be a finite Galois $2$-extension  of $\Q$ and $i \equiv 
%2 \pmod 4 $. Then the $2$-primary part of  $K_{2i-2} o_L $ vanishes exactly when  $L$
%is totally imaginary and there exists an odd prime $\ell \equiv \pm
%3 \pmod 8$  such that  $L/\Q$  is unramified outside  $\{2, \infty,
%\ell \}$ with
%$\ell \equiv \pm 3 \pmod 8$. \qed
%\end{cor}

(ii) If $i$ is odd, we discuss according to whether the number field  $L$  is complex or real. \\
(a)  If $L$ is totally complex, the vanishing of
$H^2_\et(o_L^\prime,  \Z_2(i))$ is independent of  the parity of $i$
(\cite [Lemma 2.2]{Kolster03}). Hence by the above case (i), just
studied,  we  have 

\begin{proposition}
Let   $L$ be a totally imaginary $2$-extension  of $\Q$ and $i\geq 2$ . Then the \'{e}tale cohomology $H^2_{\mbox{\'{e}t}}( o_L[1/2],
\Z_2(i))$ vanishes 
precisely when $L/\Q$
is unramified outside a set of primes $\{2, \infty, \ell \}$ with $\ell
\equiv \pm 3 \pmod 8$.
\qed
\end{proposition}

(b) If  $L$ is totally real, we assume that  $L/\Q$ is a cyclic  $2$-extension. Since  $H^2_\et (\Z[1/2],
{\Z}_{2}(i))$ is trivial  
 $$H^1_\et (\Q,  {\Z}_{2}(i))/2 \cong H^1_\et (\Z[1/2],  {\Z}/2(i)).$$
Thus 
$$D_\Q^{(i)}/\Q^{\bullet 2}=<-{\overline 1},{\overline 2}> $$
by (\ref{TatekernelforQ}). 
Hence $S\cup S_2$ is $D_\Q^{(i)}$-primitive for $(\Q, 2)$  precisely
when  $S$ is contained in $\{2, \infty, \ell_1, \ell_2 \}$  with the
Frobenius automorphisms   $\sigma_{\ell_1} (\Q(\zeta_8)/\Q)$ and
$\sigma_{\ell_2} (\Q(\zeta_8)/\Q)$ generating
$\Gal(\Q(\zeta_8)/\Q)$, or equivalently $\ell_1 \not \equiv  1 \pmod
8$,  $\ell_2 \not \equiv  1 \pmod 8$ and $\ell_1 \not \equiv \ell_2
\pmod 8$. Therefore,  Theorem \ref{motivic codescent thm}
leads to the following
%$H^2_\et(o_L^\prime,
%\Z_2(i))=0$ precisely when  $L$  is unramified outside $\{2, \ell_1,
%\ell_2 \}$  with   $\ell_1 \not \equiv  1 \pmod 8$,  $\ell_2 \not
%\equiv  1 \pmod 8$ and $\ell_1 \not \equiv \ell_2 \pmod 8$.

\begin{proposition}\label{real cyclic 2-extension} 
Let   $L$ be a totally real cyclic $2$-extension  of $\Q$ and $i\geq 2$ odd.
Then the \'{e}tale cohomology $H^2_{\mbox{\'{e}t}}( o_L[1/2],
\Z_2(i))$ vanishes exactly  when  $L$  is unramified outside a set of primes $\{2,
\ell_1, \ell_2 \}$  with   $\ell_1 \not \equiv  1 \pmod 8$,  $\ell_2
\not \equiv  1 \pmod 8$ and $\ell_1 \not \equiv \ell_2 \pmod 8$.
\qed \end{proposition}

Under hypothesis  ($\mathcal{H}_i$), the above proposition holds for general real Galois  $2$-extensions  $L$  of $\Q$. 

\subsection{}
%The deviation between $K_{2i-2} o_L \otimes \Z_2$ and
%$H^2_{\mbox{\'{e}t}}(o_L[1/2], \Z_2(i))$  (Theorem \ref{RW})
%immediately leads again to the following corollary.
In this subsection, we will apply our results to find Galois  $2$-extensions  $L$  of  $\Q$ with minimal  $2$-parts  $K_{2i-2} o_L \otimes \Z_2.$ 

(i) $2i-2 \equiv 2 \pmod 8$. 
By  Theorem \ref{RW}  and the exact sequence  (\ref{H2+H2 i even}) we have a surjective map 
$K_{2i-2} o_L\otimes \Z_2 \rightarrow ( \Z/2)^{r_1(L)}$ which becomes an isomorphism precisely when 
 $H^2_{+}(o_L[1/2], \Z_2(i))=0$. 

(ii)  $2i-2 \equiv 4 \pmod 8$. 
Theorem  \ref{RW}  and the exact sequence  (\ref{kolster2})  show that the two groups $K_{2i-2} o_L \otimes \Z_2$  and  $H^2_{+}(o_L[1/2], \Z_2(i))$  are simultaneously trivial or non trivial. 

(iii) $2i-2 \equiv 6 \pmod 8$. 
By  Theorem \ref{RW}  and the exact sequence  (\ref{H2+H2 i even}) we have an isomorphism 
$K_{2i-2} o_L \otimes \Z_2 \cong H^2_{+}(o_L[1/2], \Z_2(i))$. 

 In the above three cases we are led to study the vanishing of the positive cohomology groups  $H^2_{+}(o_L[1/2], \Z_2(i))$. 
Since    $H^2_+(\Z[1/2], {\Z}_2(i))=0$  for all integer  $i \geq 2$ (cf e.g. \cite[Proposition
2.6]{Kolster03}), Corollary \ref{going up 2,i regularity} leads immediately to the following

\begin{proposition}
Let   $L$ be a finite Galois $2$-extension  of $\Q$ and $i\geq 2$. Then the positive cohomology group $H^2_{+}(o_L[1/2], \Z_2(i))$ vanishes exactly when
$L/\Q$ is unramified outside a set of primes $\{2, \infty, \ell \}$
with $\ell \equiv \pm 3 \pmod 8$.
\qed \end{proposition}

%The deviation between $K_{2i-2} o_L \otimes \Z_2$ and
%$H^2_{\mbox{\'{e}t}}(o_L[1/2], \Z_2(i))$  (Theorem \ref{RW})
%immediately leads again to the following corollary.
%
%\begin{cor}
%Let   $L$ be a totally real finite cyclic $2$-extension  of $\Q$.
%Assume that $2i-2\equiv 0 \pmod 8$. Then  the $2$-primary part of  $K_{2i-2} o_L $ vanishes 
% precisely when
%  $L$  is unramified outside a set of primes $\{2,
%\ell_1, \ell_2 \}$  with   $\ell_1 \not \equiv  1 \pmod 8$,  $\ell_2
%\not \equiv  1 \pmod 8$ and $\ell_1 \not \equiv \ell_2 \pmod 8$.
%\end{cor}

%When  $L$ is a totally real finite Galois $2$-extension  of $\Q$ and
%$i\geq 2$  is even or $2i-2\equiv 4 \pmod 8$,  we use positive cohomology.\\

% For $i$ even, the
%genus formula for $H^2_+ (o_L^S[1/2], {\Z_2}(i)))$ gives
%$$\frac{|H^2_+(o_L[1/2], {\Z_2}(i))_{G}|}{|H^2_+(\Z[1/2], {\Z_2}(i))|}=
%\frac{|H^2(o_L[1/2], {\Z_2}(i))_{G}|}{|H^2(\Z[1/2], {\Z_2}(i))|}. $$
%Moreover, by Theorem \ref{motivic codescent thm},  ${|H^2(o_L[1/2],
%{\Z_2}(i))_{G}|}={|H^2(o_F[1/2], {\Z_2}(i))|}$ exactly when
 %$S\cup S_2$ is $D_\Q^{(i)}$-primitive
  %for $(\Q, 2)$. As seen before, this is  equivalent to the fact that $S$ is contained in a set of primes $\{2,  \ell \}$ with $\ell $
%such that  $\ell \equiv \pm 3 \pmod 8$.

%It remains to deal with the totally real case with $i\geq 2$  even.
%
%In the totally real case,  he above results are not satisfactory.  To deal with this issue, we use positive cohomology.

\begin{cor}
Let   $L$ be a totally real finite Galois $2$-extension  of $\Q$ and
  $2i-2 \not\equiv 0  \pmod 8$. Then \\
(i)  either $2i-2\equiv 2 \pmod 8$ and the surjective map $K_{2i-2} o_L\otimes \Z_2 \rightarrow ( \Z/2)^{r_1(L)}$ is an isomorphism; \\
(ii)  or  $2i-2\equiv 4 \pmod 8$ and $K_{2i-2} o_L \otimes \Z_2 =0$; \\
(iii)  or $2i-2\equiv 6 \pmod 8$ and $K_{2i-2} o_L \otimes \Z_2 =0$; \\
precisely when  $L$  is unramified outside a set of primes $\{2,  \infty, \ell \}$ with
$\ell \equiv \pm 3 \pmod 8$.
\qed \end{cor}

Finally, if $2i-2\equiv 0 \pmod 8$,  we need to discuss according to whether the Galois  $2$-extension  $L$  is real or imaginary.  
In this case, the Chern characters  
$$ch_{i,k}^\M : K_{2i-k} (o_F^S)   \longrightarrow  H^k_\M (o_F^S, {\Z}(i))$$
are isomorphisms (Theorem \ref{RW}). 

When  $L$  is complex, the vanishing of  $K_{2i-2} o_L \otimes \Z_2$  turns out to be independent of  $i \geq 2$:
\begin{cor}
Let   $L$ be a totally imaginary finite Galois $2$-extension  of $\Q$.
Then, for  $i \geq 2$,  the $2$-primary part of  $K_{2i-2} o_L $ vanishes 
precisely when  $L$  is unramified outside a set of primes $\{2,  \infty, \ell \}$ with
$\ell \equiv \pm 3 \pmod 8$.
\qed \end{cor}

When  $L$  is real, we use Proposition \ref{real cyclic 2-extension} for the vanishing of  $K_{2i-2} o_L \otimes \Z_2$. 

\begin{proposition}
Let   $L$ be a totally real finite cyclic $2$-extension  of $\Q$.
Assume that $2i-2\equiv 0 \pmod 8$. Then  the $2$-primary part of  $K_{2i-2} o_L $ vanishes 
 precisely when
  $L$  is unramified outside a set of primes $\{2,
\ell_1, \ell_2 \}$  with   $\ell_1 \not \equiv  1 \pmod 8$,  $\ell_2
\not \equiv  1 \pmod 8$ and $\ell_1 \not \equiv \ell_2 \pmod 8$.
\qed \end{proposition}
Note again that under hypothesis  ($\mathcal{H}_i$), the above proposition holds for general real Galois  $2$-extensions  $L$  of $\Q$.

%\vspace{10mm}
 \begin{tabular}{l l l l }
J. ASSIM & & & A. MOVAHHEDI \\
Universit\'e Moulay Ismail  & & & XLIM UMR 6172 CNRS/Univ. de Limoges \\
Math\'ematiques & & & Math\'ematiques et informatique\\
B.P 11201 Zitoune Mekn\`es  & & & 123, Avenue A. Thomas\\
 Mekn\`es 50000 &  & & 87060 Limoges \\
 Maroc &  & & France \\
j.assim@fs.umi.ac.ma &  & &
 mova@unilim.fr
 \end{tabular}

\end{document}